\RequirePackage{amsmath}
\RequirePackage{amssymb}
\RequirePackage{amsthm}

\newcommand{\RR}{\mathbb{R}}
\newcommand{\NN}{\mathbb{N}}

\newcommand{\ZZ}{\mathbb{Z}}
\newcommand{\CC}{\mathbb{C}}

\newtheorem{theorem}{Theorem}
\newtheorem{lemma}{Lemma}
\newtheorem{corollary}{Corollary}
\newtheorem{proposition}{Proposition}
\newtheorem{remark}{Remark}

\documentclass[dvips]{article}

\usepackage{epsfig}

\author{\textbf{V. Barutello - S. Terracini \footnote{Work partially supported by M.I.U.R. project ``Metodi Variazionali ed Equazioni Differenziali Nonlineari".}} 
\vspace{.5cm} \\
\small{ \itshape{Dipartimento di Matematica e Applicazioni}} \\
\small{ \itshape{Universit\`a di Milano Bicocca}}\\
\small{ \itshape{via Bicocca degli Arcimboldi, 8 - Milano}}\\ 
\small{ \rmfamily {e-mail:  vivina@matapp.unimib.it}}\\  
\small{ \rmfamily {\hspace{1.1cm} suster@matapp.unimib.it}}}

\title{Action minimizing orbits in the $n$-body problem with simple choreography constraint}
\date{}
\begin{document}

\maketitle

\begin{abstract}
In 1999 Chenciner and Montgomery found a remarkably simple choreographic motion for the planar 3-body problem (see \cite{CM}). In this solution 3 equal masses travel on a eight shaped planar curve; this orbit is obtained minimizing the action integral on the set of simple planar choreographies with some special symmetry constraints. In this work our aim is to study the problem of $n$ masses moving in $\RR^d$ under an attractive force generated by a potential of the kind $1/r^\alpha$, $\alpha >0$, with the only constraint to be a simple choreography: if $q_1(t),\ldots,q_n(t)$ are the $n$ orbits then we impose the existence of $x \in H^1_{2 \pi}(\RR,\RR^d)$ such that 
\begin{equation*}
q_i(t)=x(t+(i-1) \tau ),\,\,\,i=1,\ldots,n,\,\, t \in \RR,
\end{equation*}
where $\tau = 2\pi / n$. In this setting, we first prove that for every $d,n \in \NN$ and $\alpha>0$,  the lagrangian action attains its absolute minimum on the planar circle. Next we deal with the problem in a rotating frame and we show a reacher phenomenology: indeed while for some values of the angular velocity minimizers are still circles, for others  the minima of the action are not anymore rigid motions. 
\end{abstract}

{\bf Subj-class:} Dynamical Systems, Variational Methods.

{\bf MSC-class:} 70F10, 70G75, 37C80.

\section{Introduction}
In the recent years, several new periodic solutions to the $n$-body problem have been found by minimizing the action functional in space of loops having some special symmetries (see \cite{AC-Z,BC-Z,C,C2,Ch1,Ch2,CGMS,CM,CV,Ma,Mont,TV,V}). The present work is motivated by the results in \cite{FT} where the authors prove the existence of collisionless minimizers in very general classes of symmetric loops: a natural question is whether these minimizers are rigid motions corresponding to central configurations. 

To be more precise we are interested in the $2\pi$-periodic solutions to the dynamical systems of $n$ bodies, $q_1,\ldots,q_n$ with masses $m_1,\ldots,m_n$, interacting according to the Newton's law
\begin{equation}
\label{sistema dinamico}
-m_i \ddot{q}_i = \sum _{\scriptsize
\begin{array}{c}
	i,j=1 \\ i \neq j
\end{array}
}^n \nabla V_{ij}(q_i - q_j),\,\,\,\, i=1,\ldots,n
\end{equation}
where 
\begin{equation}
\label{potential}
V_{ij}(q_i - q_j) = - \frac{1}{|q_i - q_j|^\alpha}, \,\,\,\, \alpha \in \RR^*_+  
\end{equation}
is the potential that generates an attractive force between each pair of bodies.\\ 
According to \cite{AC,CGMS} a \emph{simple choreographies} is a motion in which the bodies lie on the same curve and exchange their mutual positions after a fixed time. If $q_1(t),\ldots,q_n(t)$ are the $n$ orbits then we impose the existence of $x \in H^1_{2 \pi}(\RR,\RR^d)$ such that 
\begin{equation}
\label{cor}
q_i(t)=x(t+(i-1) \tau ),\,\,\,i=1,\ldots,n,\,\, t \in \RR,
\end{equation}
where $\tau = 2\pi / n$. Moreover, we assume that our system consists of particles of equal masses. Without  loss of generality, we may assume that each $m_i=1$.

Solving (\ref{sistema dinamico}) under the constraint (\ref{cor}) is equivalent to seek the critical points of the functional
\begin{equation}
\label{azionevecchia}
{\cal A}(x) = \frac{1}{2}\sum_{h=0}^{n-1} \int_0 ^{\tau} |\dot{x}(t+h\tau)|^2 dt + \frac{1}{2}
\sum _{\scriptsize
\begin{array}{c}
	h,l=0 \\ h\neq l
\end{array}
}^{n-1}
\int_0 ^{\tau} \frac{dt}{|x(t+l\tau)-x(t+h\tau)|^\alpha}
\end{equation}
on the set 
\[
\Lambda = \{x \in H^1_{2\pi}(\RR,\RR^d) : x(t) \neq x(t+h\tau), \forall t \in \RR, h=1,\ldots,n-1 \}.  
\]

The choreography constraint implies that the quantities $\int_0 ^{\tau} \frac{dt}{|x(t+l\tau)-x(t+h\tau)|^\alpha}$ only depend on $l-h$, since there are only $n-1$ distinct average on $[0,2\pi]$ of distances between two bodies; we can then rewrite the action in (\ref{azionevecchia}) in the easier way
\begin{equation}
\label{azionenuova}
{\cal A}(x) = \frac{1}{2}\int_0 ^{2\pi} |\dot{x}(t)|^2dt +\frac{1}{2} \sum_{h=1}^{n-1}\int_0 ^{2\pi} \frac{dt}{|x(t)-x(t+h\tau)|^\alpha}.
\end{equation}

Recall that a {\it central configuration} is a configuration in which the acceleration vector of every particles is $\ddot{q}_i = -kq_i$, $k>0$, and a {\it relative equilibrium motion} is a rigid and uniform rotation of a central configuration around the center of mass; our main result is the following 
\begin{theorem}
\label{thm1}
For every $\alpha \in \RR^*_+$ and $d\geq 2$, the absolute minimum of ${\cal A}$ on $\Lambda$ is attained on a relative equilibrium motion associated to the regular $n$-gon.
\end{theorem}
In particular Theorem \ref{thm1} implies that, among simple choreographies:
\begin{itemize} 
\item [(i)] the minimum is a planar motion;
\item [(ii)] the minimum is a motion of relative equilibrium;
\item [(iii)] the minimum has minimal period $2\pi$.
\end{itemize}
In other words, we prove that the absolute minimum of the action integral among $2\pi$-periodic loops in $H^1$ is the relative equilibrium solution whose central configuration is an absolute minimum of the potential among all the configurations that generate a simple choreography. This result agrees with the one of Chenciner and Desolneux in \cite{CD} in the sense that, among choreographies, the minimum of the action is the relative equilibrium motion correspondent to the minimum among central configurations; moreover it agrees with the work of Albouy and Chenciner in \cite{AC} since the minimum is a planar motion.

In literature we know few works where explicit minima of the action are determined on special classes of loops.
From \cite{CD} we can understand the following assertion
\begin{equation}
\begin{array}{lr}
\mbox{\it If the minimal central configuration is planar and it is a regular} & \\  	
\mbox{\it polygon, then the correspondent motion of relative equilibrium } & \mbox{\it(A)}\\ 
\mbox{\it (i.e. the circle) is the absolute minimum of the action under } & \notag\\
\mbox{\it the choreography constraint.} &
\end{array}
\end{equation}
It's easy to prove that (A) is fulfilled for 3-body in any dimension. If $d \geq 3$ the minimal central configuration is never planar as soon as $n \geq 4$, therefore (A) never holds (see Correction to \cite{CD}, \cite{M}, \cite{P}). Indeed, F. Pacella in \cite{P}, proves that the planar central configurations of 4 equal masses are saddle points for the potential energy and R. Moeckel, in \cite{M}, shows that for every planar configuration of $n$ different masses, there exists a variation on the direction perpendicular to the plane of the configuration that makes the potential decrease.

Concerning the flat central configurations we know from A. Albouy, in \cite{A}, that the square is the minimal central configuration for 4 bodies and again (A) holds true. On the hands, we recall that Palmore in 1976, in \cite{Pal}, claimed that the $n$-gon is the minimum among flat central configurations, but some years later Slaminka and Woerner in \cite{SW} show the $n$-gon is not the minimum for $n \geq 6$, $n$ even. This result shows a strong obstruction to the fulfillment of (A) also in the planar case. Numerically, D. Ferrario, shows in \cite{F} that the square and the pentagon are the absolute minima for 4 and 5 bodies respectively, while if $n \geq 6$ the $n$-gon is not anymore the absolute minimizer; but if we take into account just central configurations that generates a simple choreography, we have in \cite{F} a confirmation of our result.

The above discussion shows that Theorem \ref{thm1} relies more on the choreographical assumption than on the properties of the central configurations. In order to obtain a better understanding of the structure underlying  Theorem \ref{thm1}, we proceed with the analisys of the same problem in a rotating frame with angular velocity $\omega > 0$. In this case the lagrangian action is 
$$
{\cal A}(y) = \frac{1}{2}\int_{0}^{2\pi}|\dot{y}(t)+i\omega y(t)|^2 dt + \frac{1}{2} \sum_{h=1}^{n-1}\int_0 ^{2\pi} 
\frac{dt}{|y(t)-y(t+h\tau)|^\alpha}
$$
where $y \in \Lambda$ is the loop in the rotating system on which the bodies move with the choreography constraint $y_{i+1}(t)=y_{i}(t+\tau)$. In order the deal with coercive functionals we shall choose $\omega \notin \NN$. Our results show a strong dependence on the value of the angular velocity $\omega$.
Indeed consider the function $h: \RR^*_+ \rightarrow \NN$, $h(\omega)=\min_{n \in \NN}\frac{(\omega -n)^2}{n^2}$ and let $\omega ^* := \max\{\omega : h(\omega)=1 \}$. We will prove that for $\omega \in (0,\omega^*) \backslash \{ 1 \} $ Theorem \ref{thm1} still holds true. 
In a similar manner, we shall show that for $\omega $ in appropriate $\epsilon$-neighborhood of integers $k$ coprime with $n$, $\omega \notin \NN$, the minimum of the action is still attained on the rigid equilibrium motions associated to the $n$-gon, but the minimal period of the motion of each body is now $2\pi/k$ (see Theorem \ref{int1}).

The assertion of Theorem \ref{thm1} can not be valid for every value of $\omega$: in fact if we consider values of $\omega$ close to an integer $k$ such that $g.c.d.(k,n)={\tilde k}>1$, ${\tilde k} \neq n$, the minima of the action are more complex orbits that are not anymore rigid motions, but $2\pi$-periodic curves with winding number $k$ about the center of mass; furthermore this motions show a clustering phenomenon (see Theorems \ref{kcopn} and \ref{kdivn}). Our rigorous results cover a range of $\omega$ close to the set of integers; for this reason we investigate numerically the behavior of the minimizers for values of $\omega$ close to an half integer. The numerical experiments show that the minimizers are not anymore rigid motions, moreover they are not planar provided the number of bodies is large enough ($n \geq 12$ fits) while they seem to be planar circles  when $n \leq 9$. As a final remark non planar motions appear as mountain pass points for low number of bodies ($n=3$). 
 
\vspace{.5cm}
The paper is structured as follows: the first four sections are devoted to the proof of Theorem \ref{thm1}: in Section \ref{sec2} we recall its proof in the Keplerian case, while in Sections \ref{sec3} and \ref{eigprob} there are all the elements to make the proof of Section \ref{sec5}. In Section \ref{sec7} we study the minima of the action functional when the bodies lie in a rotating system with angular velocity $\omega >0$. In the last section  we show examples for Theorems \ref{thm1}, \ref{int1}, \ref{kcopn}, \ref{kdivn} and some non-planar saddle points we have found numerically.
\section{Preliminaries}
\label{sec2}
In this section we recall, for the reader's convenience, two classical inequalities and we show the equivalent result of Theorem \ref{thm1} in the Keplerian case. This result follows from the application of Jensen Inequality as first remarked in \cite{DG} and  \cite{DGM}. 

Let $q \in H^1_{2\pi}(\RR,\RR)$  such that $[q]=\frac{1}{2\pi}\int_0 ^{2\pi} q(t)dt =0$, then
\begin{equation}
\label{poin}
\int_0 ^{2\pi} |\dot{q}(t)|^2dt \geq \int_0 ^{2\pi} |q(t)|^2dt.
\end{equation}
The inequality (\ref{poin}) is known as the \emph{Poincar\'e Inequality} and to prove it, consider 
the Fourier representation of $q$, $q(t)=\sum_{n \in \ZZ ^*} c_n e^{in t}$, $c_{-n}=\bar{c}_n$; then $\dot{q}(t)=\sum_{n \in \ZZ ^*} in c_n e^{in t}$ and 
\[
\int_0 ^{2\pi} |\dot{q}(t)|^2dt = 2 \pi \sum_{n \in \ZZ ^*} n^2 c_n ^2 \geq 2 \pi \sum_{n \in \ZZ ^*}c_n ^2 = \int_0 ^{2\pi} |q(t)|^2dt;
\]
the equality holds if and only if $c_n=0$, for every $n \notin \{-1,+1 \}$, i.e. $q(t)=c_1 e^{it} +c_{-1} e^{-it} = a \cos t + b\sin t$, $a,b \in \RR$.

If we have a vector-valued function $q \in H^1_{2\pi}(\RR,\RR^d)$, $d \in \NN \backslash \{0\}$, $q(t)=(q_1(t),\ldots,q_n(t))$, $[q]=0$, then inequality (\ref{poin}) have to be verified from every component $q_i(t)$, $i=1,\ldots,d$ and the equality is attained if and only if $q_i(t)=a_i \cos t + b_i \sin t$, $a_i,b_i \in \RR$.

Consider now a convex function $f$ on $(a,b) \subset \RR$ and $g \in L^1(0,2 \pi)$ such that $a<g(t)<b$ for every $t \in (0,2 \pi)$, then the \emph{Jensen Inequality} holds
\begin{equation}
\label{Jensen}
f \left( \frac{1}{2\pi}\int_0 ^{2\pi} g(t)dt \right) \leq \frac{1}{2\pi}\int_0 ^{2\pi} f(g(t)) dt.
\end{equation}
Observe that in equation (\ref{Jensen}) the equality holds if and only if $g(t)=cost$.

Take now an inertial system with a fixed body at the origin and a second one moving on a $2\pi$-periodic orbit $q \in H^1_{2 \pi}(\RR,\RR^d)$ according to the Newton's law
\begin{equation}
\label{sistema dinamico Kep}
-\ddot{q}_i = \nabla V(q),
\end{equation}
where $V(q) = - \frac{1}{|q|^\alpha}$, $\alpha \in \RR^*_+$, is a generalization of the Newtonian potential. 

To seek solutions of (\ref{sistema dinamico Kep}), we look for the critical points of the functional
\begin{equation}
\label{azioneK}
{\cal A}(q)= \frac{1}{2}\int_0 ^{2\pi} |\dot{q}(t)|^2dt + \int_0 ^{2\pi} \frac{dt}{|q(t)|^\alpha}
\end{equation}
defined on the set
\[
\Lambda_K = \{q \in H^1_{2\pi}(\RR,\RR^d) : [q]=0, q(t) \neq 0, \forall t \in \RR \}.  
\]
Using (\ref{poin}) and  (\ref{Jensen}) we can state that the minimun of (\ref{azioneK}) is attained on a circle. In fact, applying first the Jensen Inequality to the convex function $1/s^{\alpha /2}$, we obtain
\[
\frac{1}{2 \pi} \int_0 ^{2\pi} \frac{dt}{|q(t)|^\alpha} \geq \left( \frac{1}{\frac{1}{2 \pi} \int_0 ^{2\pi} |q(t)|^2 dt} \right)^{\alpha /2} = \left( \frac{2 \pi}{\int_0 ^{2\pi} |q(t)|^2 dt} \right)^{\alpha /2}
\]
then, using the Poincar\'e Inequality we have
\begin{eqnarray}
{\cal A}(q) & \geq & \frac{1}{2}\int_0 ^{2\pi} |\dot{q}(t)|^2dt + \frac{(2 \pi)^{\alpha/2 +1}}
                 {\left[ \int_0 ^{2\pi} |q(t)|^2 dt \right]^{\alpha/2}} \label{dis1}\\
            & \geq & \frac{1}{2}\int_0 ^{2\pi} |q(t)|^2dt + \frac{(2 \pi)^{\alpha/2 +1}}
                 {\left[ \int_0 ^{2\pi} |q(t)|^2 dt \right]^{\alpha/2}}. \label{dis2}
\end{eqnarray}
In particular 
\begin{equation}
\label{in}
{\cal A}(q) \geq \min _{s>0}\left\{ f(s) \right\}.
\end{equation}
where $f(s)= \frac{s}{2} + \frac{(2 \pi)^{\alpha/2 +1}}{s^{\alpha/2}}$, $s>0$. If there exists $\bar{q} \in \Lambda_K$ such that  in (\ref{dis1}), (\ref{dis2}), (\ref{in}) the equalities are reached, then $\bar{q}$ is the absolute minimum of the action $\cal{A}$.
Since the minimun of the function $f$ is attained at the point $s_{min}=2 \pi \alpha ^{\frac{2}{\alpha + 2}}$,
the absolute minimum of $\cal{A}$ have to verify the conditions:
\begin{equation*}
(H) \quad \left\{
\begin{array}{l}
\int_0 ^{2\pi} |\bar{q}(t)|^2dt = s_{min},\\
|\bar{q}(t)| = c \in \RR^*_+ , \\
\bar{q}_i(t) = a_i \sin t + b_i \cos t, \,\,i=1,\ldots,d. 
\end{array} \right.
\end{equation*}
If $d=3$, the only function that satisfies (H) is the planar circle $\bar{q}(t) = \alpha ^{\frac{1}{\alpha + 2}}e^{it}$. If $d=4$, we identify $\RR^4 \cong \CC^2$ and we have that (H) are satisfied if and only if $q(t)=(r_1 e^{it},r_2 e^{it})$, $r_1,r_2 \in \RR_+^*$. Since $r_2 r_1 e^{it} - r_1r_2 e^{it}=0$, we deduce that $q$ lies in $\CC$ and it is a circle. If $d>4$, we can reduce the problem to the two cases we have already examined.
\section{Some inequalities}
\label{sec3}
To simplify the notation of this section, we define the integers $K:=n-1$, $N:=K-1$.
\begin{proposition}
\label{prop1}
For every $\beta>0$, fixed any $\mu_1,\ldots \mu_K \in \RR^*_+$ the problem
\[
\min \sum_{h=1}^{K} \frac{1}{s_h^{\beta}}
\]
has a unique solution on the set $\left\{ s_h : \sum_{h=1}^{K}\mu_h s_h =1;\,\, s_1,\ldots,s_K \in \RR^*_+ \right\}$.
\end{proposition}
\begin{proof}
We want to prove that the function $f(s):=\sum_{h=1}^{K} \frac{1}{s_h^{\beta}}$ has a unique minimum under the constraint 
\begin{equation*}
\label{constraint}
\sum_{h=1}^{K}\mu_h s_h =1,
\end{equation*}
that is equivalent to prove that the function
\begin{equation*}
\label{fV}
f_C(s_1,\ldots,s_N):=\sum_{h=1}^{N} \frac{1}{s_h^{\beta}}+\mu_K ^{\beta} \left( 1- \sum_{h=1}^{N} \mu_h s_h \right)^{-\beta}  ,
\end{equation*}
has a unique minimum. This function is defined, continuous and coercive on the bounded set 
$$
{\cal S} = \left\{ (s_1,\ldots,s_N): s_i>0,i=1,\ldots,N, \sum_{i=1}^N \mu_i s_i <1  \right\}.
$$ 
Moreover $f_C$ is strictly convex on ${\cal S}$, since composition of a strictly convex function with a linear one. The uniqueness of the minimum follows. 
\\
\end{proof}

\begin{corollary}
\label{cor1}
For every $\beta>0$, fixed any $(\tilde{\mu}_1,\ldots \tilde{\mu}_K) \in (\RR^*_+)^K$ the problem
\begin{eqnarray}
&& \min \Phi (s), \quad s=(s_1,\ldots,s_K)\in (\RR^*_+)^K, \label{pcor1}\\
&& \Phi (s) := \left( \sum_{h=1}^{K} \frac{1}{s_h^{\beta}} \right) \left( \sum_{h=1}^{K} \tilde{\mu}_h s_h \right)^{\beta} \notag
\end{eqnarray}
has a unique solution, in the sense that if $m = \min \Phi(s)$ and $\Phi(\tilde{s})=m$ then $\Phi(s)=m$ if and only if there exists a positive $\lambda$ such that $s=\lambda \tilde{s}$. Moreover the minimum is attained on $\tilde{s} \in (\RR^*_+)^K$ such that 
\begin{equation}
\label{rel}
\tilde{\mu}_h = \frac{1}{\tilde{s}_h^{\beta+1}\sum_{h=1}^{K} \tilde{s}_h^{\beta}}
\end{equation}
and its value is $\sum_{h=1}^{K} \frac{1}{\tilde{s}_h^{\beta}}$.
\end{corollary}

\begin{proof}
The function $\Phi$ is homogeneous, for every $\lambda >0$ and $s \in (\RR^*_+)^K$, $\Phi(\lambda s)=\Phi(s)$. Problem (\ref{pcor1}) is equivalent to the one we have studied in Proposition \ref{prop1}, since, fixed $\tilde{\mu} \in (\RR^*_+)^K$, then $(\RR^*_+)^K = \bigcup_{\lambda >0} \lambda s^1$, where $s^1 \in (\RR^*_+)^K$
is such that $\sum_{h=1}^{K} \tilde{\mu}_h s^1_h = 1$. The unicity of the solution in the sense specified follows then from Proposition \ref{prop1}.\\ 
To find an explicit solution for (\ref{pcor1}), we impose $\frac{\partial \Phi}{\partial s_h} = 0$, $h=1,\ldots,K$ and we obtain relation (\ref{rel}).
\\
\end{proof}

\begin{proposition}
\label{dis}
For every $k \geq 2$, $k \in \NN$ and $ x \in (0,2 \pi)$  the following inequality holds
\begin{equation}
\label{trig}
1 - \cos(kx) < k^2 (1 - \cos x).
\end{equation}
\end{proposition}

\begin{proof}
Let $f(x):=1 - \cos(kx)$ and $g(x):=k^2(1 - \cos x)$, real functions defined on $(0,2\pi)$. First of all observe that the function $\varphi:=g-f$ is such that $\varphi(\pi - x)=\varphi(\pi + x)$, so we only have to prove (\ref{trig}) on the set $(0,\pi)$. Moreover, since the function  $g$ is increasing on $(0,\pi)$, $f$ increases in $(0,\pi/k)$ and $f(x)\leq 2 = f(\pi/k) < g(\pi/k)$ our claim is to show (\ref{trig}) only on the set $(0,\pi/k)$. To do that it is enought to remark that the fuction $\varphi $ is positive on the interval $(0,\pi/k)$, in fact
$$
\varphi(0)=0,\quad \varphi'(0)=0
$$
and, if $x \in (0,\pi/k)$ 
\begin{equation}
\varphi''(x)=k^2(\cos x -\cos kx)>0 \quad \mbox{for} \quad k \geq 2. 
\end{equation}
\\
\end{proof}
\section{An eigenvalue problem}
\label{eigprob}
For every $x \in \Lambda$ let $\xi^x=(\xi^x_1,\ldots,\xi^x_{n-1}) \in (\RR^*_+)^{n-1}$ defined as 
$$
\xi^x_h := \int_0 ^{2\pi} |x(t)-x(t+h\tau)|^2 dt, \,\,\,\, h=1,\ldots,{n-1},
$$
and let $\bar{\xi}_h := \xi^{\bar{x}}_h = 8 \pi \sin^2(\pi h/n)$, where $\bar{x}$ is the circle of radius one centered at the origin.\\ 
The positive numbers  $\bar{\mu}_h$ are defined as 
\[
\bar{\mu}_h = \frac{1}{c\bar{\xi}_h ^{\alpha/2 + 1}}, \,\,\, \text{where} \,\,\, c= \sum_1 ^{n-1} \frac{1}{\bar{\xi}_h ^{\alpha/2}}.
\]
In what follows $x_i$, $i=0,\ldots,n-1$, are $n$ vectors of $\RR^d$. If $j \geq n$ then $x_j=x_i$, where $i=j \mbox{ mod } n$, so when we write $x_j$ with $j \geq n$, we actually want to write the correspondent $x_i$.

Define the linear operator $\Delta^\mu:  \RR^{nd} \rightarrow \RR^{nd}$ in the following way
\begin{equation}
\label{delta}
(\Delta^\mu x)_i := \sum_{h=1}^{n-1} \bar \mu_h (2x_i - x_{i+h} - x_{i-h}),
\end{equation}
where $x_i \in \RR^d$, $i=0,\ldots,{n-1}$. \\
In the first part of this section our aim is to study the eigenvalues and the eigenvectors of $\Delta^\mu$. Remark that $\Delta^\mu$ acts on each component of every $x_i \in \RR^d$ in the same way
\begin{equation}
\label{deltacomp}
(\Delta^\mu x)_i^j := \sum_{h=1}^{n-1} \bar \mu_h (2x_i^j - x_{i+h}^j - x_{i-h}^j)
\end{equation}
for $i=0,\ldots,{n-1}$ and $j=1,\ldots,d$. For this reason we can associate to $\Delta^\mu$ the square matrix $nd \times nd$
\begin{equation*}
\left[
\begin{array}{ccc}
D^1 & & 0 \\
    & \ddots	& \\
0 & & D^1    
\end{array}
\right]
\end{equation*}
where $D^1$ is the $n \times n$ matrix associated to the operator $\Delta^\mu$ restricted to the first components. The eigenvalues of $\Delta^\mu$ are the eigenvalues of the matrix $D^1$ and their multiplicity is $d$ times the multiplicity of the ones of $D^1$. If $v_l \in \RR^n$ is an eigenvector of $D^1$ correspondent to the eigenvalue $\delta_l$ then the operator $\Delta^\mu$ has an eigenspace correspondent to $\delta_l$ spanned by the vectors 
$$ 
\left( v_l,0_{\RR^n},\ldots,0_{\RR^n} \right), \left( 0_{\RR^n},v_l,0_{\RR^n},\ldots,0_{\RR^n} \right),\ldots,
\left( 0_{\RR^n},\ldots,0_{\RR^n}, v_l \right) \in \RR^{nd}
$$
Using these remarks we will have a complete description of $\Delta^{\mu}$ by studying the eigenvalues and the eigenvectors of the matrix $D^1$.

\begin{lemma}
\label{studiodidelta}
The matrix $D^1$ defined below has $[n/2]+1$ distinct eigenvalues
\begin{equation}
\label{eigenvalues}
\delta_0=0,\,\,\,\delta_1=\frac{1}{2\pi},\,\,\,\delta_l<\frac{l^2}{2\pi}\,\,\,l=2,\ldots,[n/2].
\end{equation}
The multiplicity of $\delta_0$ is 1 and the correspondent eigenspace is the one spanned by a constant vector of $\RR^n$.\\
The multiplicity of $\delta_1,\ldots,\delta_{[n/2]}$ is 2 except when $n$ is even, in this case the multiplicity of $\delta_{n/2}$ is 1.\\
The eigenspace correspondent to the eigenvalue $\delta_l$, $l=1,\ldots,[n/2]$, is spanned by the pair of vectors of $\RR^n$
\begin{equation}
\label{eigenvectors}
\left(\sin \left( \frac{2 \pi}{n}il \right)\right)_{i=0,\ldots,n-1}, \,\, \left(\cos \left( \frac{2 \pi}{n}il \right)\right)_{i=0,\ldots,n-1}.
\end{equation}
\end{lemma}

\begin{remark}
\label{even}
If $n$ is even, the multiplicity of $\delta_{n/2}$ is 1 since  $\left(\sin \left( \frac{2 \pi}{n}il \right)\right)_i=(\sin(\pi i) )_i=(0)_i$.
\end{remark}

\begin{proof}
$\delta_0=0$ and its multiplicity is 1 since the kernel of $D^1$ is not trivial and it is composed by the vectors $\lambda {\bf 1}_{\RR^n}$, where ${\bf 1}_{\RR^n}=(1,\ldots,1)$ and $\lambda \in \RR$.

The vectors in (\ref{eigenvectors}), that we denote with $v_l$ are eigenvectors for $D^1$, in fact for every $l=1,\ldots,[n/2]$
$$
(D^1 v_l)_i = \left[2\sum_{h=1}^{n-1} \bar \mu_h (1- \cos \left( \frac{2 \pi}{n}hl \right) )\right] (v_l)_i
$$
and then 
\begin{equation}
\label{deltal}
\delta_l =  2\sum_{h=1}^{n-1} \bar \mu_h (1- \cos \left( \frac{2 \pi}{n}hl \right)), \,\,\, l=1,\ldots[n/2].
\end{equation}
We can calculate $\delta_1$:
\begin{eqnarray*}
2 \sum_{h=1}^{n-1} \bar \mu _h (1 - \cos \left( \frac{2 \pi}{n}h \right)) &=& 2 \sum_{h=1}^{{n-1}} \frac{2\sin^2(\pi h /n)}{c\bar{\xi}_h ^{\beta + 1}}= \frac{4}{c} \sum_{h=1}^{{n-1}} \frac{\bar{\xi}_h}{8 \pi \bar{\xi}_h ^{\beta + 1}}\\
&=&  \frac{1}{2 \pi c} \sum_{h=1}^{{n-1}} \frac{1}{\bar{\xi}_h ^\beta}= \frac{1}{2 \pi};
\end{eqnarray*}
if we substitute the inequality of Proposition \ref{dis} in (\ref{deltal}) we obtain the estimate
\begin{equation}
\label{stima}
\delta_l < \frac{l^2}{2\pi},
\end{equation}
for $l=2,\ldots,[n/2]$.

The eigenspaces corresponding to $\delta_1, \ldots ,\delta_{[n/2]}$ are 2-dimensional (except the case $n$ even, where the eigenspace correspondent to  $\delta_{[n/2]}$ is 1-dimensional, see Remark \ref{even}) since 
\begin{equation}
\label{orto}
\sum_{i=1}^{n} \sin \left( \frac{2 \pi}{n}il \right) \cos \left( \frac{2 \pi}{n}il \right) = 
\frac{1}{2} \sum_{i=1}^{n} \sin \left( 2\frac{2 \pi}{n}il \right)=0;
\end{equation}
some therms in the second sum of (\ref{orto}) are zero ($\sin (k\pi)$), the others are pairs of  opposite numbers of the kind 
$$
\sin \left( 2l \frac{2 \pi}{n} i \right) = - \sin \left( 2l \frac{2 \pi}{n} (n-i) \right).
$$
We conclude the proof, since we have found a set of $n$ independent eigenvectors for the $n$-dimensional matrix  $D^1$. 
\\
\end{proof}
\begin{remark}
\label{complessi}
Each eigenspace of the operator $D^1$ (except the kernel) is spanned by the pair of vectors in (\ref{eigenvectors}). These 2-dimensional subspaces of $\RR^n$ can be interpretate as 1-dimensional subspaces of $\CC^n$ generated by the vectors
$$
\left( e^{\frac{2\pi}{n}Jil} \right)_{i=0,\ldots,n-1}, \quad l=1,\ldots,[n/2], 
$$
where $J$ is the complex unit.
When $n$ if even and the multiplicity of $\delta_{n/2}$ is 1, $v_{n/2} = \left( e^{\pi Ji}\right)_{i=0,\ldots,n-1} = (1,-1,\ldots,1,-1)$ is a real element of $\CC^n$.
\end{remark}
\begin{proposition}
\label{proppp}
Consider the periodic problem 
\begin{equation}
\label{(P)}
\left\{
\begin{array}{l}
	-\ddot{x}^j(t)=\lambda (\Delta^{\mu}x(t))^j, \quad j=1,\ldots,d\\	
	x^j \in H^1_{2\pi}(\RR,\RR^n)
\end{array}
\right.
\end{equation}
where $x(t)=\left( x^1(t),\ldots, x^d(t)\right)$, $\lambda \in \RR$ and $\Delta^\mu:  \RR^{nd} \rightarrow \RR^{nd}$ is defined in (\ref{delta}). (\ref{(P)}) admits solutions if and only if $\lambda=\lambda_{kl} = \frac{k^2}{\delta_l}$, where $k \in \ZZ$ and $\delta_l$ is an eigenvalue for $\Delta^{\mu}$. Moreover a solution of (\ref{(P)}) with $\lambda = \lambda_{kl}$ has the form
\begin{equation}
\label{sol}
x^j(t)=v_l\gamma(t), \quad j=1,\ldots,d
\end{equation}
where $v_l \in \RR^n$ is an eigenvector of $D^1$ corresponding to the eigenvalue $\delta_l$ and the scalar function $\gamma$ solves the periodic linear problem
\begin{equation}
\label{gamma}
\left\{
\begin{array}{l}
	-\ddot{\gamma}(t)=k^2\gamma(t), \\	
	\gamma(t+2 \pi) = \gamma (t), \quad t \in \RR.
\end{array}
\right.
\end{equation}
\end{proposition}
\begin{proof}
Since the operator $\Delta^{\mu}$ acts on each component in the same way, we can rewrite the first $d$ equations of 
(\ref{(P)})
\begin{equation}
\label{(Pr)}
	-\ddot{x}^j(t)=\lambda D^1 x^j(t), \quad j=1,\ldots,d.
\end{equation}
Let $N$ a $n \times n$ matrix such that $N^{-1} D^1 N = D$ is diagonal (the column of $N$ are the eigenvectors of $D^1$), then (\ref{(Pr)}) is equivalent to
\begin{equation}
\label{(Prr)}
	-\ddot{y}^j_l(t)=\lambda \delta_l y^j_l(t), \quad j=1,\ldots,d, \quad l=1,\ldots,n, 
\end{equation}
where $y^j(t) = N^{-1}x^j(t)$. Problem (\ref{(Prr)}), under the periodic boundary condition, admits solution if and only if $\lambda \delta_l = k^2$, $k \in \ZZ$; so we deduce that  to satisfy (\ref{(Prr)}) $y^{j}(t)=(0,\ldots,0,\gamma(t),0,\ldots,0)$, where $\gamma(t)$ is at the $l$-th place and verify (\ref{gamma}) with $\lambda \delta_l = k^2$. So we deduce that 
$$
x^{j}(t) = Ny^{j}(t) = v_l\gamma(t),
$$  
where $v_l$ is the eigenvector of $D^1$ correspondent to $\delta_l$ and $\gamma(t)$ solves (\ref{gamma}).
\\
\end{proof}
Let $x \in \Lambda$ and define 
\begin{equation}
\label{corpi}
x_i: \RR \rightarrow \RR^d, \quad x_i(t)=x(t+i\tau), \quad i=0,\ldots,n-1;
\end{equation}
the motion of the $(i+1)$-th body.
For $j=1,\ldots,d$, $x^j_i:\RR \rightarrow \RR$ is the $j$-th component of the function $x_i$; moreover we define
\begin{equation}
\label{xj}
x^j: \RR \rightarrow \RR^n, \quad x^j(t)=(x^j_i(t))_{i=0,\ldots,n-1}, \quad j=1,\ldots,d.
\end{equation}
that represents the $n$ $j$-th components of the functions $x_i$. Since $x \in \Lambda$, $x^j \in H^1_{2\pi}(\RR,\RR^n)$.
\begin{proposition}
\label{propultima}
Take $x \in \Lambda$ and $x^j$ as in (\ref{xj}). Then the problem (\ref{(P)}) admits a solution if and only if $\lambda = \frac{k^2}{\delta_l}$, where $\delta_l$ is an eigenvalue for $\Delta^{\mu}$ and $k=l+rn$, $r \in \ZZ$.
\end{proposition}
\begin{proof}
We impose to the general solution (\ref{sol}) the condition of simple choreography, i.e.
\begin{equation}
\label{core}
x_{i+1}(t) = x_{i}(t+\tau), \quad i=0,\ldots,n-1,
\end{equation}
where we understand $x_n(t)=x_0(t)$.\\
First of all observe that the orbit $x(t)$ lies in an even dimensional space; in fact the scalar equation (\ref{sol}) does not admit a solution under the constraint (\ref{core}).\\
If $d = 2m$, using Remark \ref{complessi}, we write the $2m$ scalar equations (\ref{sol}) as $m$ complex equations
all equal to
\begin{equation}
\label{comeq}
x_i(t) = e^{\tau J il} e^{Jkt}, \quad i=0,\ldots,n-1
\end{equation}
and the simple choreography constraint is equivalent to $e^{\tau Jl} = e^{\tau Jk}$; since $k \geq l$ we obtain $k=l+rn$, $r \in \ZZ$.\\
While $m$ equal parametric equations in $\CC^m$ determine a 2-dimensional subspace of $\CC^m$, $m$ equations equal to (\ref{comeq}) assure that the orbit $x(t)$ is planar circle and each body run along it $l$ times in a period $2\pi$.
\\
\end{proof} 
 
\begin{lemma}
\label{da dim}
For every $x \in \Lambda$
\begin{equation}
\label{da fare}
\frac{1}{2}\int_0 ^{2\pi} |\dot{x}(t)|^2dt \geq 2 \pi n \sum_{h=1}^{n-1} {\bar \mu}_h \xi^x_h
\end{equation}
and the equality is attained if and only if $x$ is a circle.
\end{lemma}

\begin{remark}
Problem  (\ref{da fare}) is an eigenvalue problem and the set of solutions is a linear space. So, if the equality in  (\ref{da fare}) is attained on a circle then it is attained on every one.
\end{remark}

\begin{proof}
First of all observe that if $x \in \Lambda$ is not the zero function, (\ref{da fare}) is equivalent to
\[
\frac{\frac{1}{2}\int_0 ^{2\pi} |\dot{x}(t)|^2dt}{ \left\langle \Delta^{\mu} x,x \right\rangle_{L^2}} \geq 2 \pi,
\]
where to apply the operator $\Delta^{\mu}$ to $x$ we consider $x(t) = (x^j(t))^{j=1,\ldots,d}$ and the scalar product in $L^2$ is $\left\langle \Delta^{\mu} x,x \right\rangle_{L^2} = \int^{2\pi}_0 \left\langle \Delta^{\mu} x,x \right\rangle_{\RR^{nd}} = \sum_{i=0}^{n-1} \int^{2\pi}_0 \left\langle \Delta^{\mu} x_i,x_i \right\rangle_{\RR^{d}}$. In fact 
\begin{eqnarray*}
\int_0 ^{2\pi} \left\langle \Delta^{\mu} x,x \right\rangle_{\RR^{nd}} &=&  \sum_{i=0}^{n-1} \sum_{h=1}^{n-1} \bar \mu_h \int_0 ^{2\pi} (2x_i^2 - x_i x_{i+h} - x_i x_{i-h})\\
&=&  \sum_{i=0}^{n-1} \sum_{h=1}^{n-1} \bar \mu_h \int_0 ^{2\pi} (2x_i^2 - 2x_i x_{i+h})\\
&=&  \sum_{i=0}^{n-1} \sum_{h=1}^{n-1} \bar \mu_h \int_0 ^{2\pi} | x_i - x_{i+h}|^2 = n \sum_{h=1}^{n-1} \bar \mu_h \xi_h^x.
\end{eqnarray*}
The functional 
\begin{equation*}
\label{J}
J(x):=\frac{\frac{1}{2}\int_0 ^{2\pi} |\dot{x}(t)|^2dt}{ \left\langle \Delta^{\mu} x,x \right\rangle_{L^2}},
\end{equation*}
defined on $\Lambda$, is bounded from below (by a strictly positive constant) and, if $(x_n)_n$ is a minimizing sequence then\\
-  $\| x_n \|_{\Lambda} = \int_0 ^{2\pi} |\dot{x_n}(t)|^2 dt$ is bounded, \\
- the compactness in $L^2$ of the space of $2\pi$-periodic functions with zero mean value in $H^1(0,2\pi)$ ensures that $x_n \rightarrow x$ in $L^2$ and then $\left\langle \Delta^{\mu} x_n, x_n \right\rangle_{L^2} \rightarrow \left\langle \Delta^{\mu} x,x \right\rangle_{L^2}$,\\ 
- the boundedness of $\| x_n \|_{\Lambda} = \| \dot{x}_n \|_{L^2}$ ensures that $\dot{x}_n \rightharpoonup \dot{x}$ in $L^2$ and then (by the lower-semi continuity of the norm)
\[
J(x) \leq {\lim \inf}_n J(x_n).
\]
So we can state that the minimun of $J$ exists. The aim is to prove that this minimum is $2\pi$.\\
If $x$ is a $2\pi$-periodic solution of 
\begin{equation*}
\label{pm}
\min_{x \in \Lambda x \neq 0}J(x),
\end{equation*}
then $(x^j)_{j=1,\ldots,d}$, defined in (\ref{xj}), satisfies problem (\ref{(P)}), where $\lambda$ assume the minimal value that assures a solution. Thanks to Proposition \ref{propultima}, the possible values for $\lambda$ are 
$$
\lambda_{l}=\frac{l^2}{\delta_l}, \,\,\,l=1,\ldots,[n/2].
$$
Using (\ref{eigenvalues}) we can conclude that 
\begin{equation*}
2\pi = \min_{x \in \Lambda x \neq 0}J(x),
\end{equation*}
and that the minimizer is the circle with minimal period $2\pi$.
\\
\end{proof}
\section{Proof of Theorem \ref{thm1}}
\label{sec5}
The idea of the proof is to find a functional $\bar{\cal A}$ defined on the set $\Lambda$ such that:
\begin{itemize}
\item[(i)] ${\cal A}(x) \geq \bar{\cal A}(x)$ for every $x \in \Lambda$;
\item[(ii)] attains its minimum if and only if $x$ is a circle;
\item[(iii)] ${\cal A}(x) = \bar{\cal A}(x)$ if and only if $x$ is a circle.
\end{itemize}
With this aim we remember the definitions of $\xi^x$, $\bar{\xi}$, $\bar{\mu}$ and $c$ given in Section \ref{eigprob} and we use the Jensen Inequality, applied to the convex function $f(s)=\frac{1}{s^{\alpha/2}}$ and $g(t)=|x(t)-x(t+h\tau)|^2$, to obtain
\begin{equation}
\label{primadis}
{\cal A}(x) \geq  \frac{1}{2}\int_0 ^{2\pi} |\dot{x}(t)|^2dt +\frac{1}{2}\sum_{h=1}^{n-1} \frac{(2\pi)^{\frac{\alpha}{2}+1}}{ (\xi^x_h)^ {\alpha/2}}.
\end{equation}
Remark that in (\ref{primadis}) the equality is attained if $x$ is a circle.\\
Take now 
\begin{equation} 
\label{Atilde}
\tilde {{\cal A}}(x) := \frac{1}{2}\int_0 ^{2\pi} |\dot{x}(t)|^2dt +  \tilde {c} \frac{1}{ \left( \sum_{h=1}^{n-1} {\bar \mu}_h \xi^x_h \right)^ {\alpha/2}}
\end{equation} 
where $\tilde {c}=\frac{(2\pi)^{\frac{\alpha}{2}+1}}{2}c$. Corollary \ref{cor1} ensures the inequality
\begin{equation}
\label{secondadis}
{\cal A}(x) \geq \tilde {\cal A} (x).
\end{equation}
The quantities $\bar{\mu}_h$ are defined in such a way that the minimum of the functional $\left( \sum_{h=1}^{n-1} \frac{1}{(\xi^x_h)^{\alpha/2}}\right) \left( \sum_{h=1}^{n-1} {\bar \mu}_h \xi^x_h \right)^ {\alpha/2}$ is obtained if $\xi^x_h = \bar{\xi}_h$ and its value is $c$, so if $x$ is a circle the equality in (\ref{secondadis}) is attained. 

From Lemma \ref{da dim} follows that if we define
\begin{equation}
\label{Abar}
\bar{\cal A}(x) := 2 \pi n \sum_{h=1}^{n-1} {\bar \mu}_h \xi^x_h + \tilde {c} \frac{1}{ \left( \sum_{h=1}^{n-1} {\bar \mu}_h \xi^x_h \right)^ {\alpha/2}}, \,\,\,\,\,\, x \in \Lambda,
\end{equation}
then 
\begin{equation}
\label{terzadis}
\tilde {\cal A}(x) \geq \bar {\cal A} (x),
\end{equation}
and the equality holds if and only if $x$ is a $2\pi$-periodic circle.\\ 
The function $\bar {\cal A}$ has a unique minimum that it is reached when $x$ is a circle. In fact let us consider the function $g:\RR^*_+ \rightarrow \RR$, $g(y)= 2 \pi n y + \frac{\tilde c}{y^{\alpha/2}}$, which has a unique minimum at $\bar y =\left(\frac{\alpha \tilde c}{4 \pi n}\right)^{\frac{1}{\frac{\alpha}{2} + 1}}$. If $y(x) = \sum_{h=1}^{n-1} {\bar \mu}_h \xi^x_h$ then $\bar {\cal A} (x)=g(y(x))$ and the minimum is attained if and only if  $x=\bar x _R$, where $\bar x _R$ is the circle of radius 
\begin{equation}
\label{radius}
R = \left(\frac{\alpha \tilde c}{4\pi n}\right)^{\frac{1}{\alpha + 2}}
\end{equation} 
This concludes the proof.
\section{Minima in rotating systems}
\label{sec7}
In this section we try to generalize the previous results to the case when  the $n$ particles move in a non-inertial system, more precisely in a system rotating with angular velocity of intensity $\omega$. We shall prove that Theorem \ref{thm1} still holds in some range of $\omega$. This is not a technical obstruction: indeed the assertion of Theorem \ref{thm1} does not hold true for many angular velocities $\omega$. 

Suppose $d=2$, let $x(t)$, the orbit on which the bodies move in the inertial system and $y(t)$ the orbit in the non-inertial one. We suppose that the two systems are linked by a rotation with angular velocity of intensity $\omega>0$, then the relation between the orbit in the inertial system and the one in the the rotating system is
$$
x(t)=e^{J\omega t}y(t).
$$
We want to minimize the action among $2\pi$-periodic loops in the rotating system, $y \in \Lambda$, which satisfies the choreography constraint $y_{i+1}(t)=y_{i}(t+\tau)$, where $y_i(t)=y(t+i\tau)$, $i=0,\ldots,n-1$, is the motion of the $(i+1)$-th body and $y_n(t)=y_0(t)$. In the inertial system the motion verifies the following conditions 
\begin{equation}
\label{c1}
x(t+2\pi)=e^{2 \pi J\omega}x(t), \quad x_{i+1}(t)=e^{-J\omega \tau}x_{i}(t+\tau).
\end{equation}

We seek the loop that minimize the action integral in the rotating system. The mutual distances among bodies are invariant under a rotation, and, since the potential depends only on these quantities, when we write the action as a function of the orbit in the rotating frame, $y(t)$, we have the only difference in the kinetic part
\begin{equation}
\label{action_rot}
{\cal A}(y) = \frac{1}{2}\int_{0}^{2\pi}|\dot{y}(t)+J\omega y(t)|^2 dt + \frac{1}{2} \sum_{h=1}^{n-1}\int_0 ^{2\pi} \frac{dt}{|y(t)-y(t+h\tau)|^\alpha}
\end{equation}

\begin{proposition}
\label{omega=kcop}
If $\omega = k \in  [2,n-1]\cap\NN$ and if $k,n$ are coprime, then the action does not admit a minimum. 
\end{proposition}
\begin{proof}
Consider the loops $y_\nu(t)=R_\nu e^{Jkt}$, $R_\nu \rightarrow +\infty$ if $\nu \rightarrow +\infty$; $(y_\nu(t))_\nu$ form a not converging minimizing sequence for ${\cal A}$.
\\
\end{proof}

\begin{proposition}
\label{omega=n}
If $\omega = n$, then the action admits an infinity of minima that form a continuum.
\end{proposition}
\begin{proof}
From (\ref{c1}), we have that the orbit in the inertial system is $2\pi$-periodic and verifies the choreography condition $x_{i+1}(t)=x_{i}(t+\tau)$. Since the center of mass is not fixed, all $2\pi$-periodic circle with radius expressed in (\ref{radius}) are minima in the inertial system. In the rotating system the minima are $y(t)=(R(\alpha,n)e^{\pm Jt}+c)e^{-Jnkt}$, $c \in \CC$.
\\
\end{proof}

\begin{proposition}
\label{omega=k}
If $\omega = k \in [2,n-1]\cap\NN$ and $g.c.d.(k,n)= \tilde{k} >1$, $\tilde{k} \neq n$, let $j:=n/\tilde{k}$ and $\tilde{j}:=k/\tilde{k}$. Then
\begin{itemize}
\item[(i)] the action does not achieve its infimum;
\item[(ii)] any minimizing sequence $(y^{(\nu)})_\nu$ has the form
\begin{eqnarray*}
&& y_i^{(\nu)} = e^{-2\pi J \frac{i}{j}}\left( c_m^{(\nu)} + \eta(t +i\tau) \right) e^{Jkt} + o(1),, \quad i=0,\ldots,n-1\\
&& i \equiv m \mod j \quad \mbox{ and }  \quad m=0,\ldots,j-1
\end{eqnarray*}
where $|c_m^{(\nu)}| \rightarrow + \infty$ if  $\nu \rightarrow + \infty$ for $m=0,\ldots,j-1$ and the function $\eta$ is the minimum of the action described in Theorem \ref{thm1} with $n=\tilde{k}$.
\end{itemize}
\end{proposition}

\begin{proof}
Conditions  (\ref{c1})  in the inertial system imply that
\begin{equation}
\label{c2}
x(t+2\pi)=x(t), \quad
x_{i+1}(t)=e^{-2\pi J \tilde{j}/j}x_{i}(t+\tau);
\end{equation}
then we have $j$ equal systems (up to rotations) with $\tilde{k}$ bodies each in which the bodies are linked by a choreography constraint. 
We can split the action as the sum of two part 
$$
{\cal A} = j{\cal A}_{\tilde k} + {\cal A}_{int} 
$$
where ${\cal A}_{\tilde k}$ is the action of each subsystem with $\tilde{k}$ bodies and ${\cal A}_{int}$ is the sum of the interactions between each pair of bodies in different subsystems. Therefore the infimum of the action is $j$-times the action of the minimal motion with $\tilde{k}$ bodies.

On every minimizing sequence ${\cal A}_{int} \rightarrow 0$, so that the center of mass of every subsystem 
will go to infinity in a radial direction and the angle between two of such directions is a multiple of $2\pi/j$; moreover the motion on every subsystem tends to the  minimal orbit with $\tilde{k}$ bodies we have found in Theorem \ref{thm1}.
\\
\end{proof}

Consider the function $h: \RR^*_+ \rightarrow \NN$, $h(\omega)=\min_{n \in \NN}\frac{(\omega -n)^2}{n^2}$ and let $\omega ^* := \max\{\omega : h(\omega)=1 \}$, ($\omega ^* \cong 1,3333$). We can state the following 

\begin{theorem}
\label{int1}
If $\omega \in (0,\omega^*) \backslash \{+1\}$, then the functional in (\ref{action_rot}) attains its minimum on a circle with minimal period $2\pi$ and radius depending on $n$, $\alpha$ and $\omega$.
\end{theorem}

\begin{proof}
Following the proof of Theorem \ref{thm1}, we would like to show the existence of a constant depending on $\omega$, $c(\omega)$, such that the following inequality holds
\begin{equation}
\label{disrot}
\frac{1}{2}\int_{0}^{2\pi}|\dot{y}(t)+J\omega y(t)|^2 dt \geq c(\omega)n\sum_{h=1}^{n-1} {\bar \mu}_h \xi^y_h,
\end{equation}
for every loop $y \in \Lambda$. Using the same arguments of Lemma \ref{da dim} we have that (\ref{disrot}), is equivalent to
\begin{equation}
\label{disrot2}
\frac{\frac{1}{2}\int_{0}^{2\pi}|\dot{y}(t)+J\omega y(t)|^2 dt}{\left\langle \Delta^\mu y,y \right\rangle_{L^2}} \geq c(\omega)
\end{equation}
and that if $y$ is a $2\pi$-periodic solution of
$$ 
\min_{y \in \Lambda\,\,y\neq 0}\frac{\frac{1}{2}\int_{0}^{2\pi}|\dot{y}(t)+J\omega y(t)|^2 dt}{\left\langle \Delta^\mu y,y \right\rangle_{L^2}} 
$$
then $y$ satisfies the linear system
\begin{equation}
\label{(PR)}
\left\{
\begin{array}{l}
	-\ddot{y}(t) +2J\omega\dot{y}+\omega^2y=\lambda (\Delta^{\mu}y(t)),\\	
	y \in H^1_{2\pi}(\RR,\CC)
\end{array}
\right.
\end{equation}
where $\lambda \in \RR$ and $\Delta^\mu:  \RR^{2n} \rightarrow \RR^{2n}$ is defined in (\ref{delta}). Following Proposition \ref{proppp}, the periodic problem (\ref{(PR)}) admits solutions if and only if $\lambda = \lambda_{kl} = \frac{\mu_k^\omega}{\delta_l}$, where $\delta_l$ is an eigenvalue for $\Delta^\mu$ and $\mu_k^\omega = (\omega - k)^2$, $k \in \ZZ$, is an eigenvalue for the linear second order periodic problem in $\CC$
\begin{equation}
\label{sop}
\left\{
\begin{array}{l}
	\ddot{\gamma}(t) + 2J\omega \dot{\gamma}(t) + \omega^2 \gamma(t)=\mu \gamma(t),\\	
	\gamma(t+2\pi) = {\gamma}(t), \quad t \in \RR.
\end{array}
\right.
\end{equation}
A solution of (\ref{(PR)}) with $\lambda = \lambda_{kl}$ has the form $ y(t)=v_l\gamma_k(t)$, where $v_l \in \CC$ is an eigenvector of $D^1$ corresponding to the eigenvalue $\delta_l$ and the complex function $\gamma_k$ solves (\ref{sop}) with $\mu = \mu_k^\omega$, i.e. $\gamma_k (t) = e^{ikt}$.\\
The choreography condition in the rotating system, $y_{i+1}(t) = y_i(t+\tau)$ implies that the relation between $k$ and $l$ is $k = l + rn$, $r \in \ZZ$, so that the admissible eigenvalues for (\ref{(PR)}) are $\lambda_l=\frac{(\omega-(l + rn))^2}{\delta_l}$. \\
The constant $c(\omega)$ is the minimim value of $\lambda_l$. Since we have for which taken $\omega \in (0,\omega^*) \backslash \{+1\}$ and $\delta_l \leq l^2/2 \pi$, then the minimal value for $\lambda_l$ is verified for $l=1$ we obtain 
$$
c(\omega)=2\pi (\omega -1)^2.
$$
We can then conclude that the minimal loop  is a $2\pi$-periodic circle with radius $R_\omega = \left(\frac{\alpha \tilde{c}}{2 n c(\omega)}\right)^{1/(\alpha +2)}$.
\\
\end{proof}

\begin{remark}
\label{ddim}
If the bodies move in a $d$-dimensional space, we can use the same arguments of Proposition \ref{propultima} to show that the orbit is planar and, precisely, that it lies on the plane orthogonal to the direction of the rotation. In fact, we can decompose the motion in $[d/2]$ 2-dimensional spaces (if $d$ is odd, there is no motion in the last dimension since the choreography constraint could not be verified) such that the first of them is the plane of the rotation. The eigenvalues correspondent to the motions in the 2-dimensional subspaces are all $2\pi$ except the one of the plane ortoghonal to the rotation that is $c(\omega)$; since $c(\omega)<2\pi$ the motion is a planar circle on the rotating plane.
\end{remark}

\begin{theorem}
\label{kcopn}
Suppose that $n$ and ${k}$ are coprime. Then there exist $\epsilon = \epsilon(\alpha,n,{k} )$ such that if $\omega \in ({k} - \epsilon, {k} + \epsilon)$ the minimum of the action is attained on a circle with minimal period $2\pi/{k}$ that lies in the rotating plane with radius depending on $n$, $\alpha$ and $\omega$.
\end{theorem}

\begin{proof}
Take $\omega \in (\omega^*,n - 1/2)$, $\omega \notin \ZZ+1/2$ and  let $[\omega]_n = {k} \in [2,n-1] \cap \NN$,
the integer closest to $\omega$. Suppose  ${k}$ coprime with $n$. We define a new linear operator $\Delta ^{\mu,{k}}$ that on the two components of the rotating plane acts in the following way
\begin{equation}
\label{deltak}
(\Delta^{\mu,{k}} y)^j_i := \sum_{h=1}^{n-1} \bar {\mu}_h^{k} (2y_i^j - y_{i+h}^j - y_{i-h}^j),
\end{equation}
where  $\bar {\mu}_h^{k}= \frac{1}{c^{k}(\bar{\xi}^{k}_h) ^{\alpha/2 + 1}}$,  $\bar{\xi}^{k}_h = 8 \pi \sin^2(\pi {k}h/n)$ and $c^{k}= \sum_1 ^{n-1} \frac{1}{(\bar{\xi}^{k}_h) ^{\alpha/2}}$.
For the other components $\Delta^{\mu,{k}}$ acts in the same way as $\Delta ^{\mu}$. The eigenvalues for the operator $\Delta ^{\mu,{k}}$ corresponding to the components on the rotating frame are the same eigenvalues of the matrix $D^1$, but rearranged in a different way, we will call them $\delta^{{k}}_l$. It is easy to calculate that $\delta ^{k}_{k}=1/2\pi$. To prove that the minimum of the action is a circle with minimal period $2\pi/k$ in the rotating plane we would like to verify that
\begin{equation}
\label{min}
\min_l \frac{(\omega-l)^2}{\delta ^{ k}_l} = \frac{(\omega-{ k})^2}{\delta ^{ k}_{ k}} = 2\pi(\omega-{ k})^2.
\end{equation}
Since the minimal eigenvalue is $1/2\pi$, then condition (\ref{min}) is implied by 
\begin{equation}
\label{min2}
\delta^{ k}_{\max} \leq \frac{1}{2\pi(\omega-{ k})^2} \min_k \left((\omega-{ k}+1)^2,(\omega-{ k}-1)^2\right),
\end{equation}
where $\delta^{ k}_{\max}$ is the largest eigenvalue of the linear operator $\Delta^{\mu,{k}}$. Fixed $n$ and $\alpha$, (\ref{min2}) is verified if we take $\omega$ in an appropriate $\epsilon$-neighborhood of ${ k}$. To conclude the proof remark that, when condition (\ref{min}) is verified, then, with the same argument we used in Remark \ref{ddim}, we can assert that the minimum is the planar circle with minimal period $2\pi/{ k}$ that lies in the rotating plane. Moreover the radius of a minimal circle is $R_\omega = \left(\frac{\alpha \tilde{c}}{2 n c(\omega)}\right)^{1/(\alpha +2)}$, where $c(\omega)=2\pi(\omega-{ k})^2$.
\\
\end{proof}

\begin{remark}
\label{3corpi}
If $n=3$, $\delta^{2}_{\max} = 1/2\pi$ and condition (\ref{min2}) is always verified. In this particular case if $\omega \in \{ 1,2 \}$, the action does not admit a minimum; if $\omega = 3/2$, the minima are all circles with minimal period $\pi$ or $2\pi$; if $\omega \in (0,3/2) \backslash  \{ 1 \}$ the minimum is attained on circles with period $2\pi$ and if $\omega \in (3/2,3) \backslash  \{ 2 \}$ the minimum is on circles with period $\pi$.
\end{remark}

\begin{remark}
\label{deltamax}
Calculating the values of $\delta^{ k}_{\max}$, we remark that if $4 \leq n \leq 9$, condition (\ref{min2}) is verified for every $\alpha >0$, this implies that if we take ${ k}$ coprime with $n$, then for every $\omega \in ({ k} - 1/2, { k} + 1/2)$ the minimum of the action is attained on a circle with minimal period $2\pi/{ k}$ that lies in the rotating plane. 
\end{remark}

\begin{theorem}
\label{kdivn}
Take ${k} \in \NN$ and $g.c.d.(k,n)= \tilde{k} >1$, $\tilde{k} \neq n$. Then there exists $\epsilon = \epsilon(\alpha,n,{k})>0$ such that if $\omega \in ({k} - \epsilon, {k} + \epsilon) \backslash \{k \}$ the minimum of the action is attained on a planar $2\pi$-periodic  orbit with winding number ${k}$ which is not a relative equilibrium motion. 
\end{theorem}

\begin{proof}
Let $(\omega_\nu)_\nu$, $\omega_\nu \rightarrow {k}$ as $\nu \rightarrow + \infty$. If $m(\omega_\nu)$ is the minimal value of the action correspondent to $\omega_\nu$ then  $m(\omega_\nu) \rightarrow \inf j{\cal A}_{\tilde k}$ (see Proposition \ref{omega=k}) for the continuity of the action in $\omega$. Then there exists $\nu_0$ such that if $\nu > \nu_0$, $m(\omega_\nu)$ is attained on a orbit of the kind described in {\em (ii)} of Proposition \ref{omega=k}.
\\
\end{proof}

\begin{remark}
\label{circan}
With the same argument we used in the proof of Theorem \ref{kdivn}, we can assert that there exists $\epsilon = \epsilon(\alpha,n)>0$ such that if $\omega \in (n - \epsilon, n + \epsilon) \backslash \{n \}$ then the minimum of the action is attained on a orbit of the kind  $y(t)=(R(\alpha,n)e^{\pm Jt}+c)e^{-Jnkt}$, $c \in \CC$ (see Proposition \ref{omega=n}).
\end{remark}

In the theorems above we have examined the minimizers of the action functional when the angular velocity of the rotating system is in the interval $(0,n-1+\epsilon)$, where $\epsilon=\epsilon(\alpha,n,n-1)$ according with Proposition \ref{omega=kcop}, \ref{omega=k} and Theorems  \ref{int1}, \ref{kcopn}, \ref{kdivn}, or when $\omega \in (n - \xi, n + \xi)$, $\xi=\xi(\alpha,n)$ according with Proposition \ref{omega=n} and Remark \ref{circan}. Consider now $\omega > n$, such that
\[
\omega = \bar \omega + ln, \quad \quad l \in \NN^*, \quad \bar \omega \in (0,n).
\]
If $\bar \omega=0$ then we are in the same situation described in Proposition \ref{omega=n}. 
If $\bar \omega = k$,  Propositions \ref{omega=kcop} and \ref{omega=k} still hold true substituting $\omega$ with $\bar \omega$. If  $\bar \omega$ is not an integer, then Theorems \ref{int1}, \ref{kcopn}, \ref{kdivn} still hold 
substituting $\omega$ with $\bar \omega$ and with the following differences: in Theorems \ref{int1} and \ref{kcopn} the minimal period of the minimizing circles are in this case respectively  $2\pi/ln$ and $2\pi/(ln+k)$; in Theorem \ref{kdivn} the winding number of the minimizing orbit is $k+ln$.

\section{Numerical results and further comments}
\label{sec6}
We give now some numerical examples of the planar orbits we described in Theorem \ref{kdivn}. The algorithm we use to minimize the action functional acts following the steepest descendent flow, for this reason it could happen that the orbits in Figures 1, 2 are just local minima for the action and not the global ones; anyway their qualitative behavior is the one described in Theorem \ref{kdivn}. In Figure 1 there are some examples of the orbits we have described in Theorem \ref{kdivn} when the angular velocity $\omega$ is closed to an integer that divides the number of bodies. In Figure 2 the we consider values of the angular velocity that are close to an integer that is not coprime with $n$ and that does not divide it.

\begin{figure}[h!]
\begin{center}
\begin{tabular}{cc}
{\psfig{figure=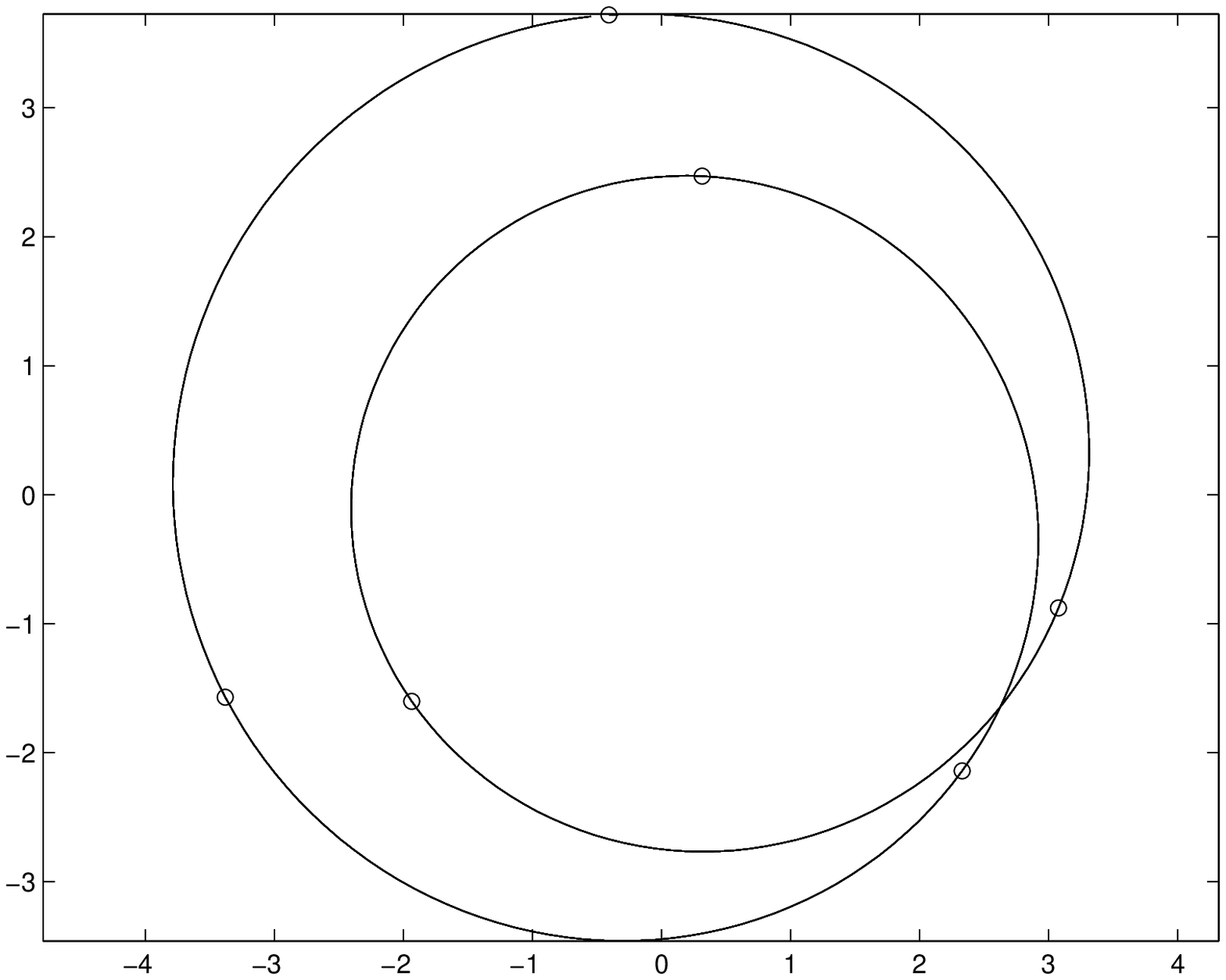,width=5.0cm}} \quad & \quad
{\psfig{figure=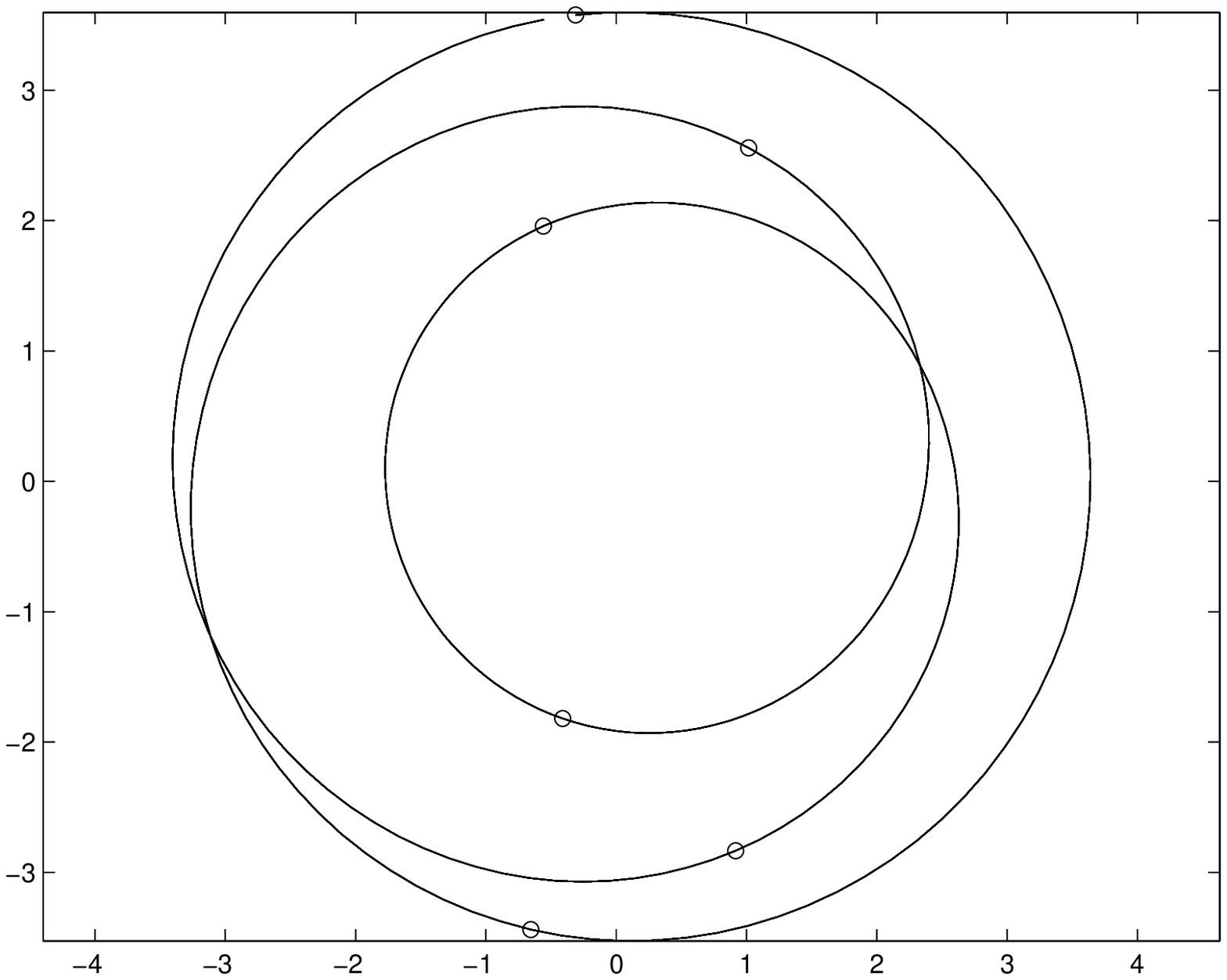,width=5.0cm}} \\
{$n=6$, $\omega = 1.8$, $\alpha=1$} & {$n=6$, $\omega = 3.2$, $\alpha=1$} \\
{\psfig{figure=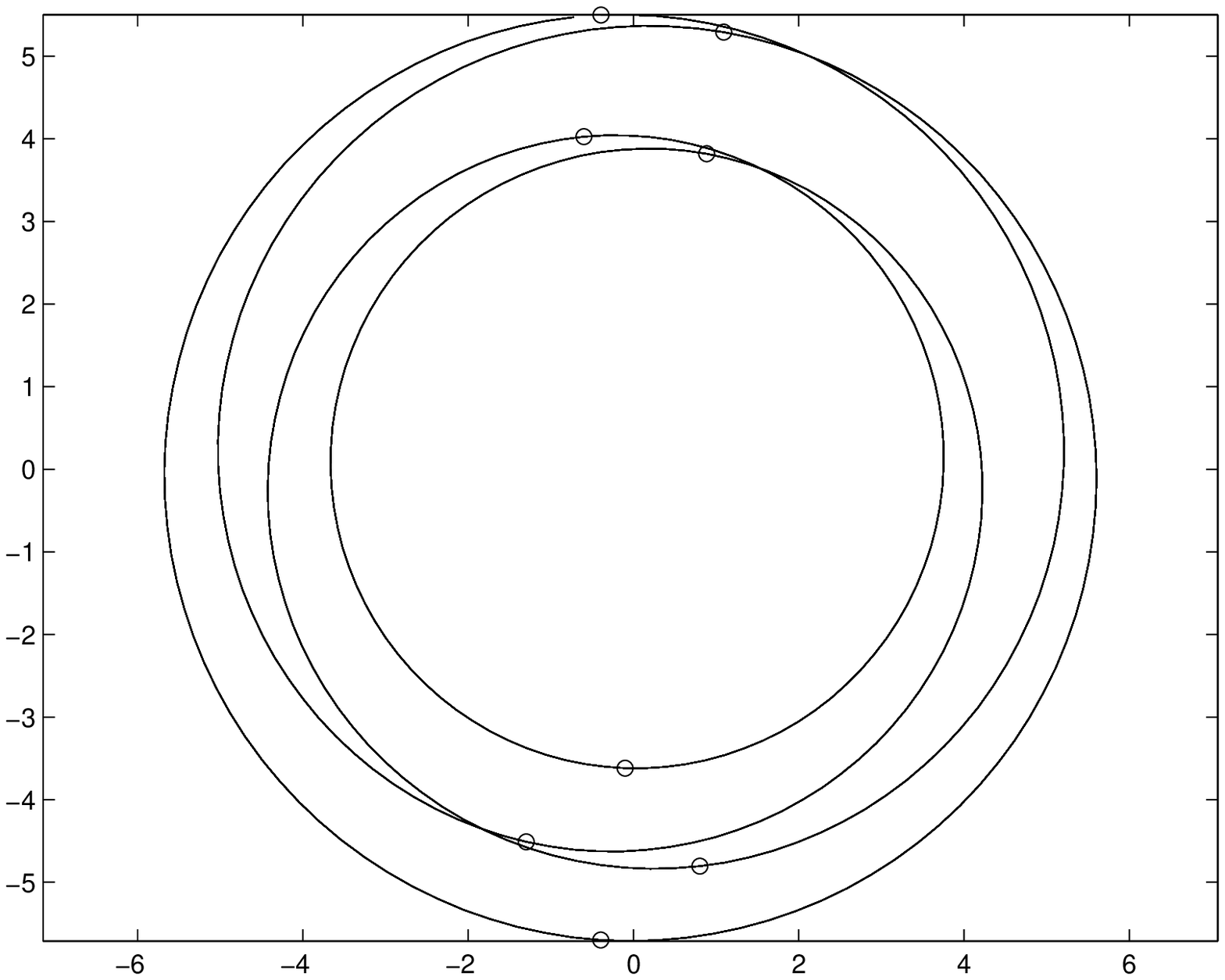,width=5.0cm}} \quad & \quad
{\psfig{figure=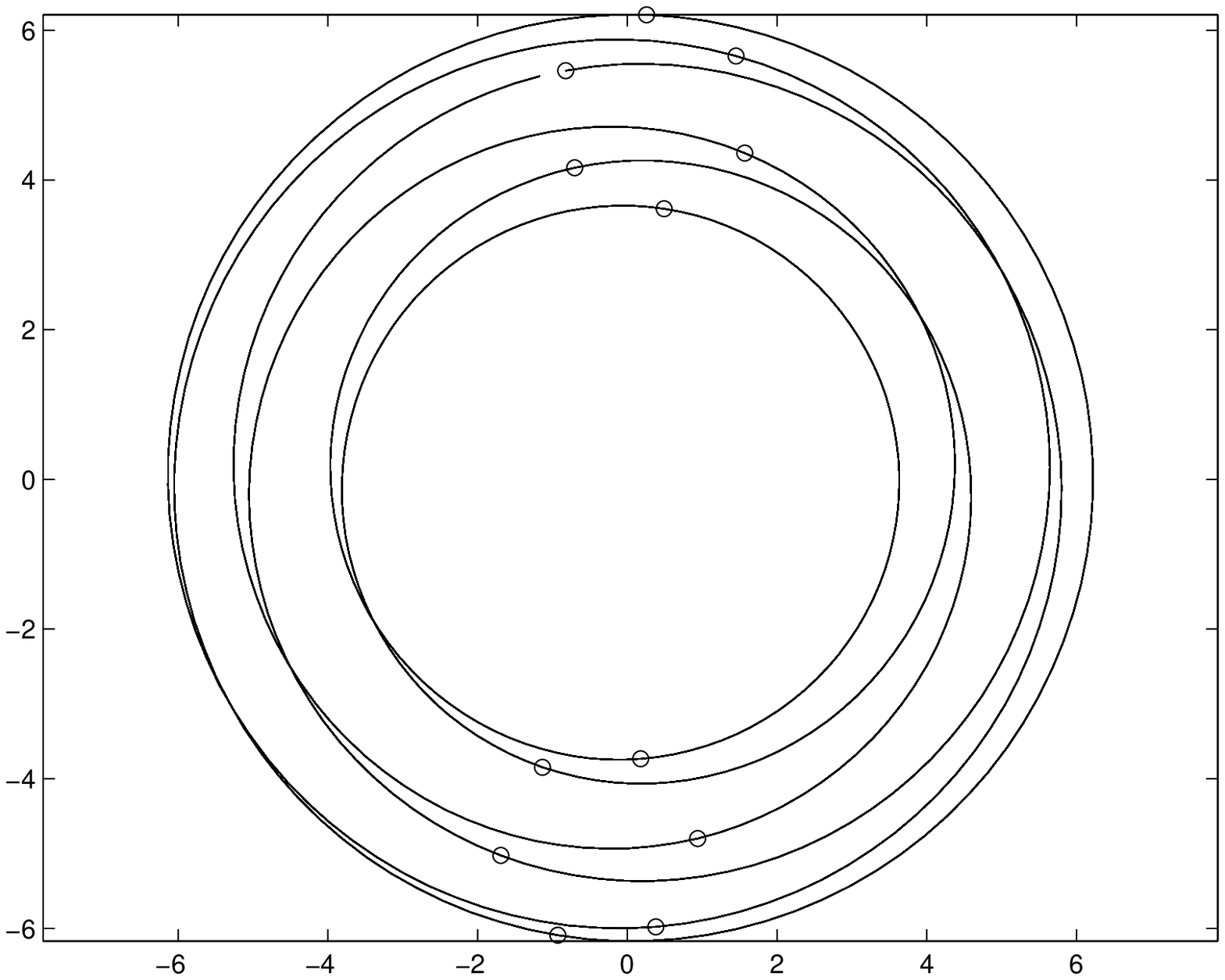,width=5.0cm}} \\
{$n=8$, $\omega = 3.9$, $\alpha=1$} & {$n=12$, $\omega = 6.1$, $\alpha=1$} 
\end{tabular}
\begin{tabular}{c}
{\psfig{figure=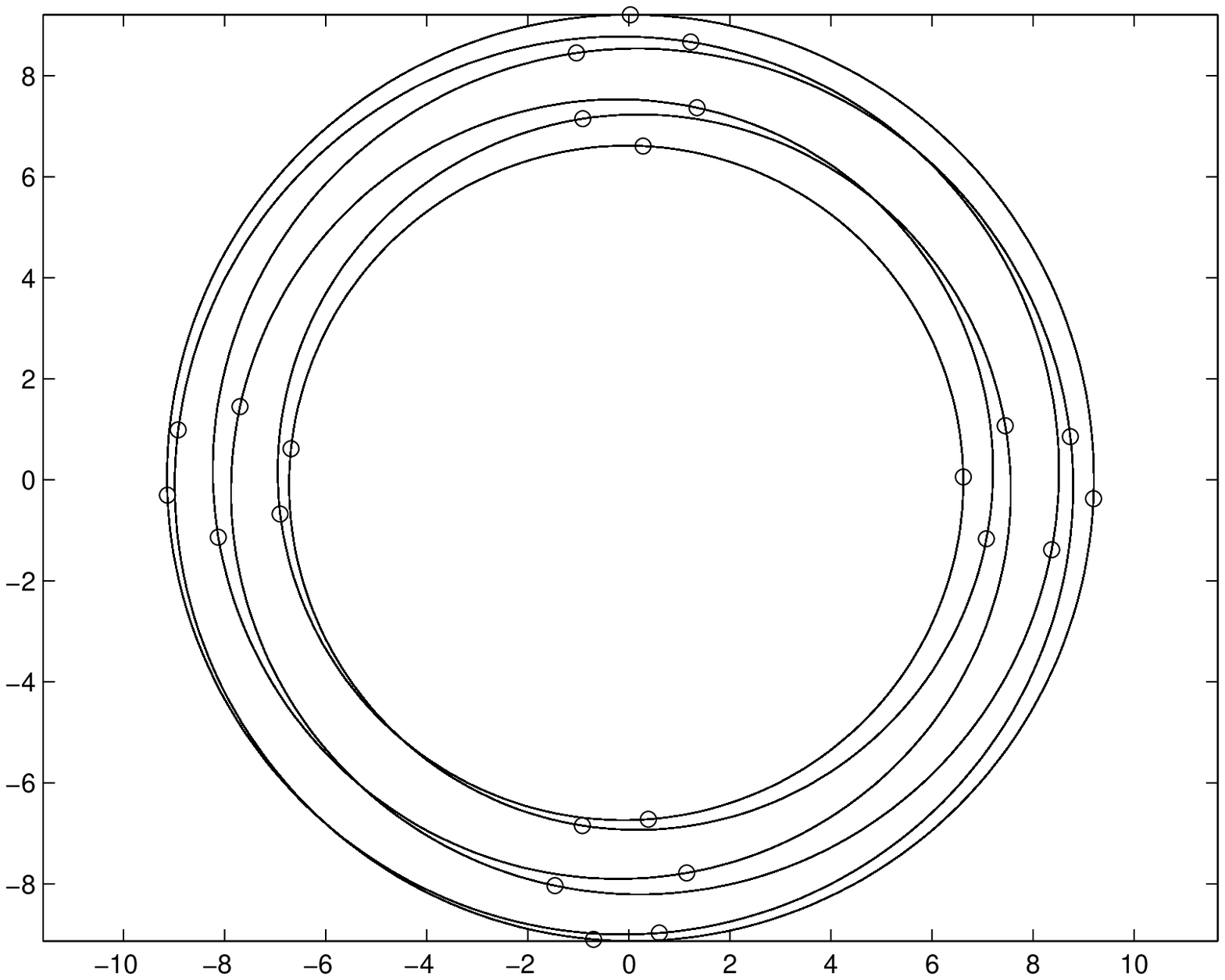,width=5.0cm}}\\
$n=24$, $\omega = 6.1$, $\alpha=1$
\end{tabular}
\caption{Examples for Theorem \ref{kdivn}, $\omega$ close to an integer that divides $n$.}
\end{center}
\end{figure}

\newpage

\begin{figure}[ht]
\begin{center}
\begin{tabular}{cc}
{\psfig{figure=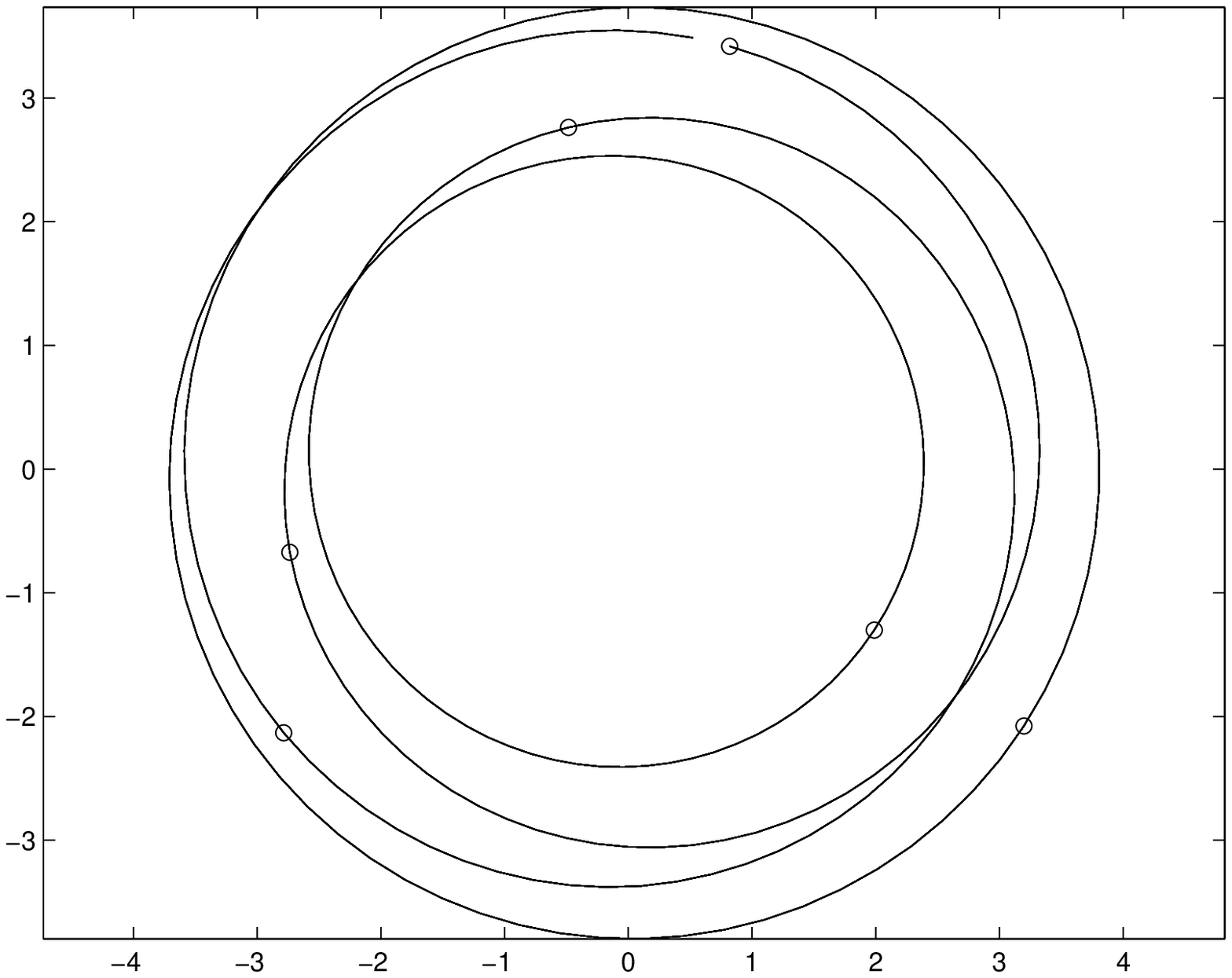,width=5.0cm}} \quad & \quad
{\psfig{figure=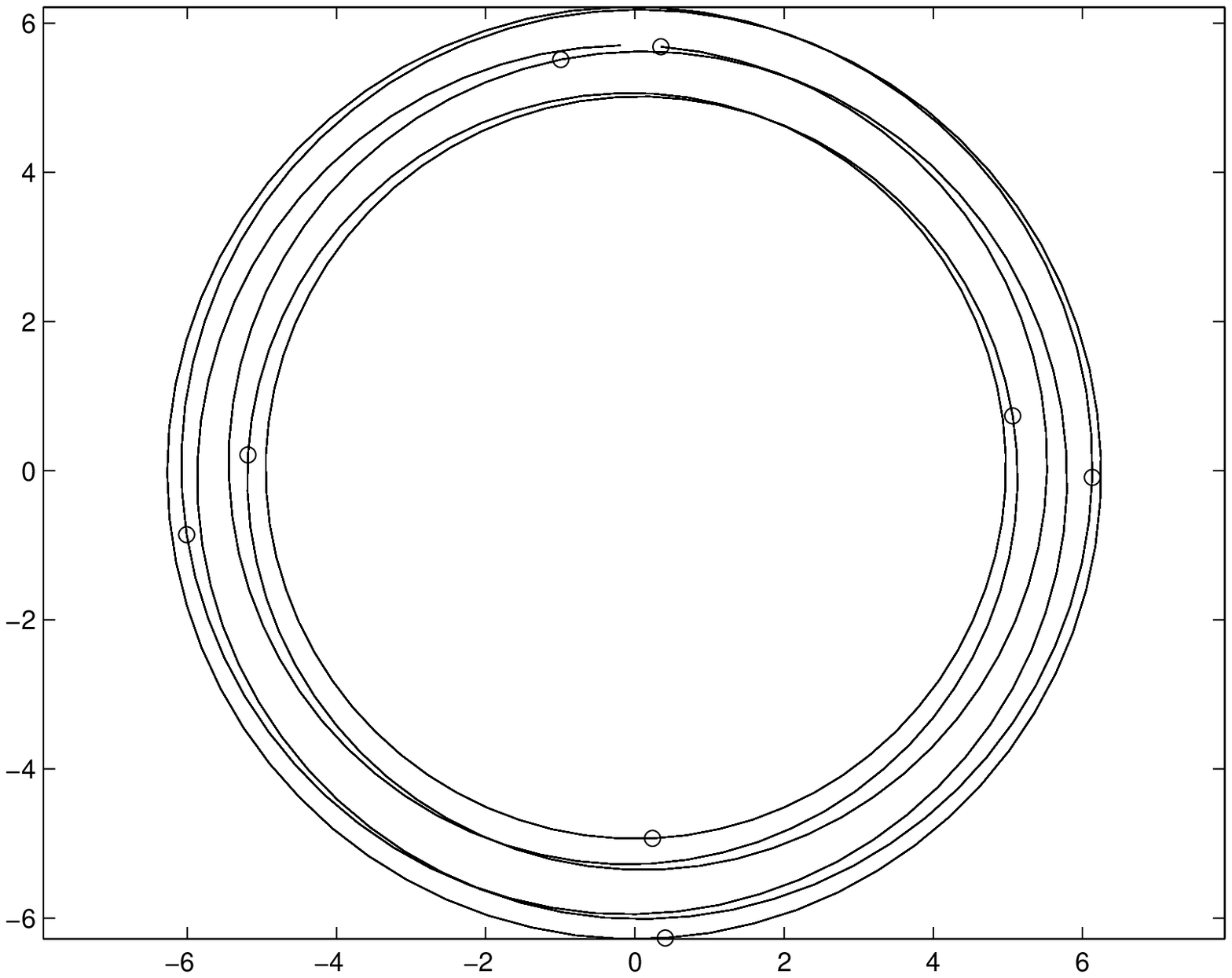,width=5.0cm}} \\
{$n=6$, $\omega = 4.2$, $\alpha=1$} & {$n=8$, $\omega = 6.1$, $\alpha=1$} 
\end{tabular}
\begin{tabular}{c}
{\psfig{figure=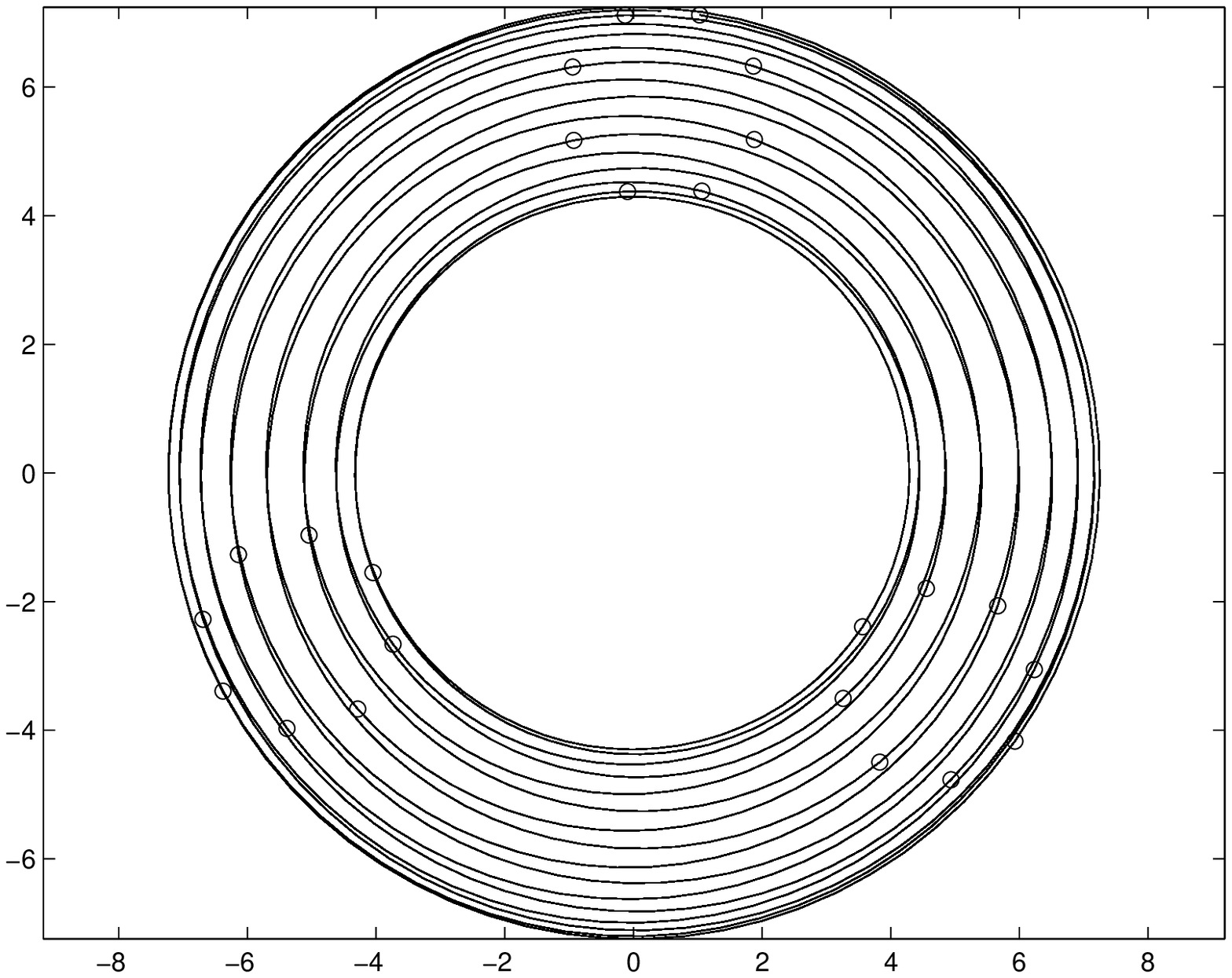,width=5.0cm}}\\
$n=24$, $\omega = 16.1$, $\alpha=1$
\end{tabular}
\caption{Examples for Theorem \ref{kdivn}, $\omega$ close to an integer that divides $n$.}
\end{center}
\end{figure}

Theorems \ref{int1}, \ref{kcopn} and \ref{kdivn} examine the minima of the action integral for values of $\omega$ close to integers; this limitation is not just a technical obstruction, we can in fact find numerical examples of minima that are not the ones of Theorems \ref{int1}, \ref{kcopn} and \ref{kdivn} when $\omega$ is close to the half of an integer. In the left column of Figure 3 we show some examples of non planar minima correspondent to values of $\omega$ close to the half on an integer; on the right side we show their projections on the plane $x_3 = 0$. Remark that in the example with 12 bodies, the projection of the curve on the horizontal plane is still a circle with winding number 7, that is the integer  closest to $\omega$; this orbit has a similar behavior to the ``hip-hop" found by A. Chenciner and A. Venturelli in \cite{CV}. In the other two pictures, with 15 and 17 bodies respectively, the projections on the $x_3=0$ plane are orbits similar to the one described on line by Sim\`o with winding number equal to the integer closest to $\omega$.

\newpage 

\begin{figure}[h!]
\begin{center}
  \begin{tabular}{cc}
  {\psfig{figure=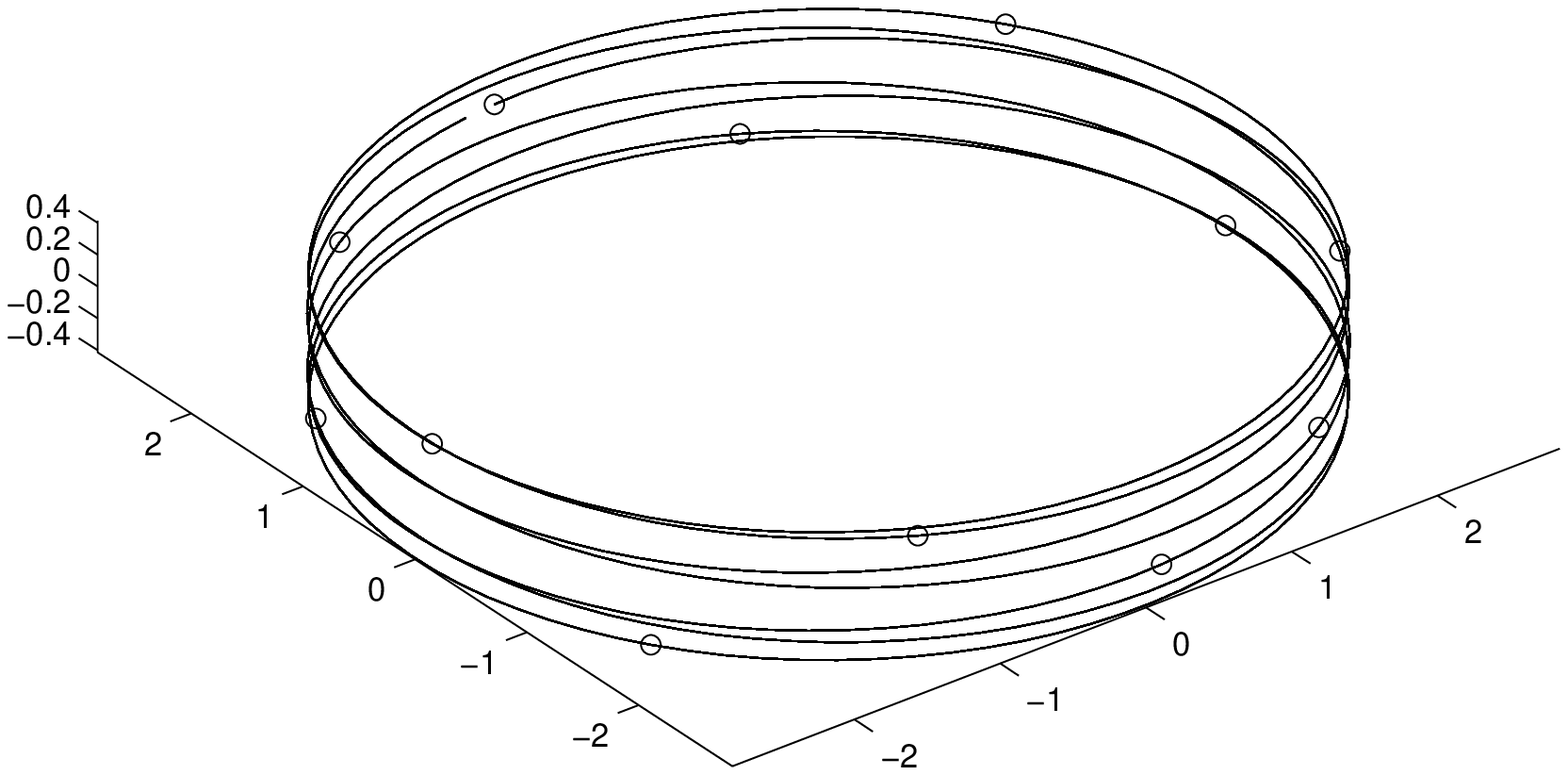,width=5.0cm}} \quad & \quad
  {\psfig{figure=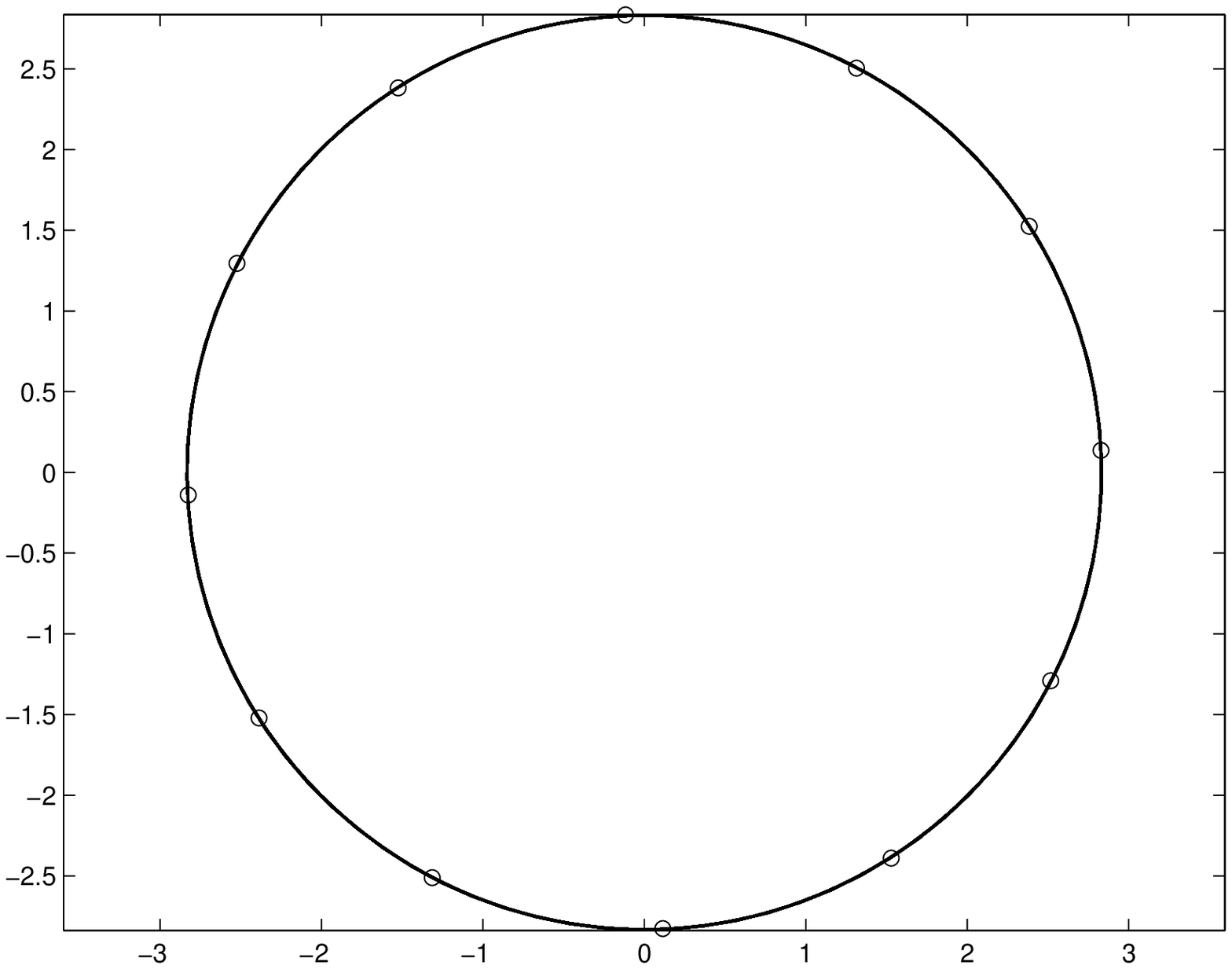,width=5.0cm}} \\
  {$n=12$, $\omega = 6.55$, $\alpha=1$} & \\
  {\psfig{figure=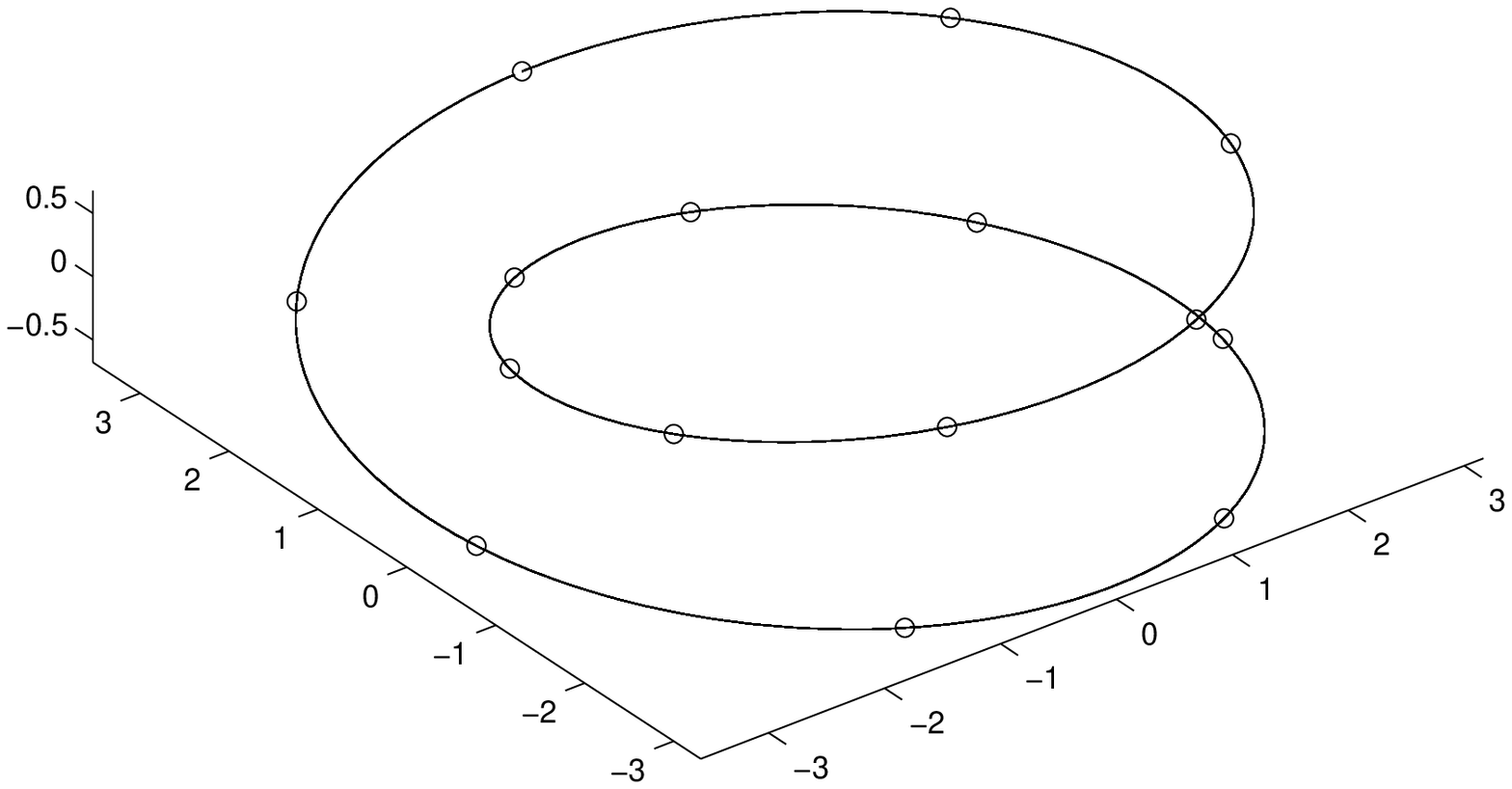,width=5.0cm}} \quad & \quad
  {\psfig{figure=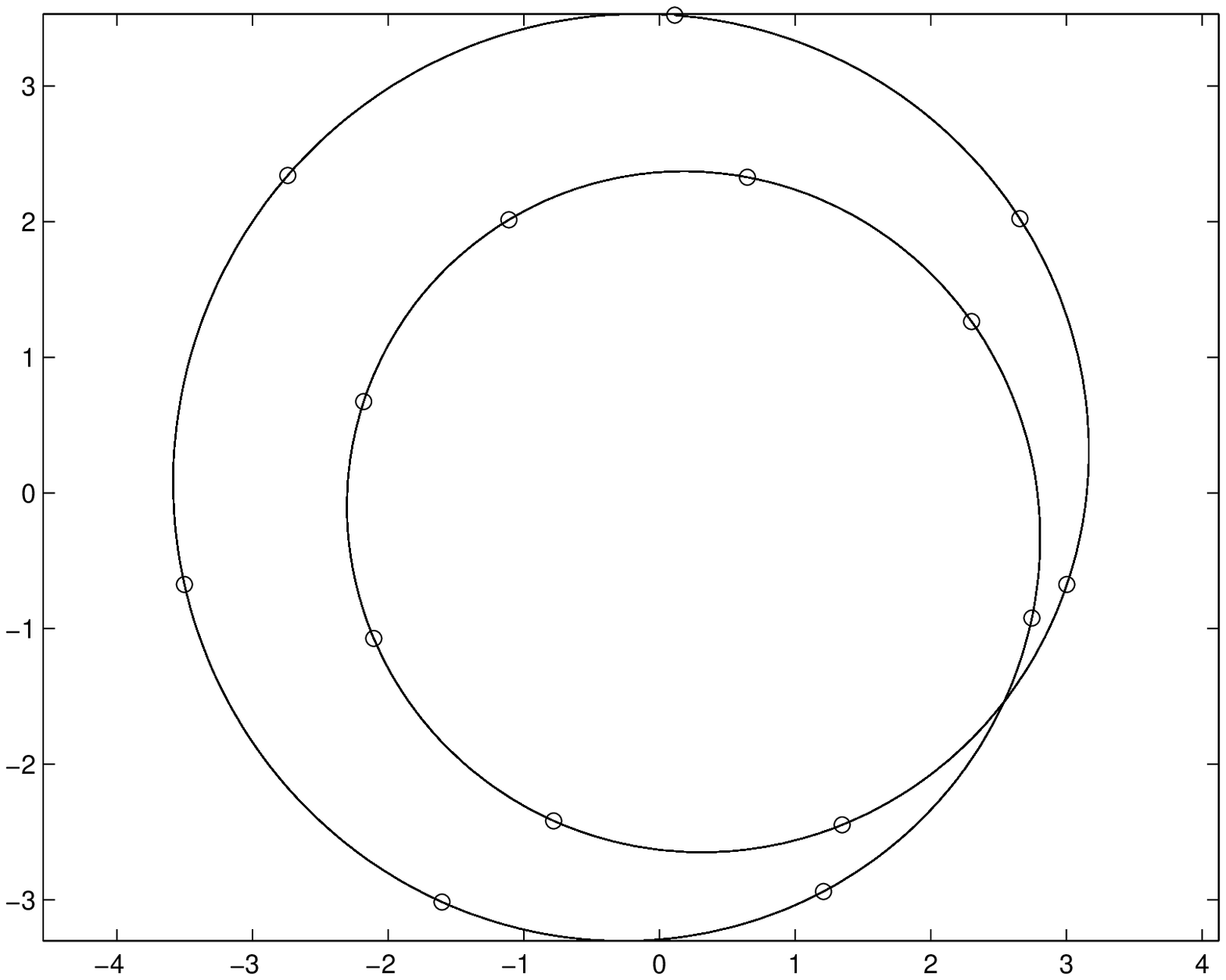,width=5.0cm}} \\
  {$n=15$, $\omega = 2.45$, $\alpha=1$} & \\
  {\psfig{figure=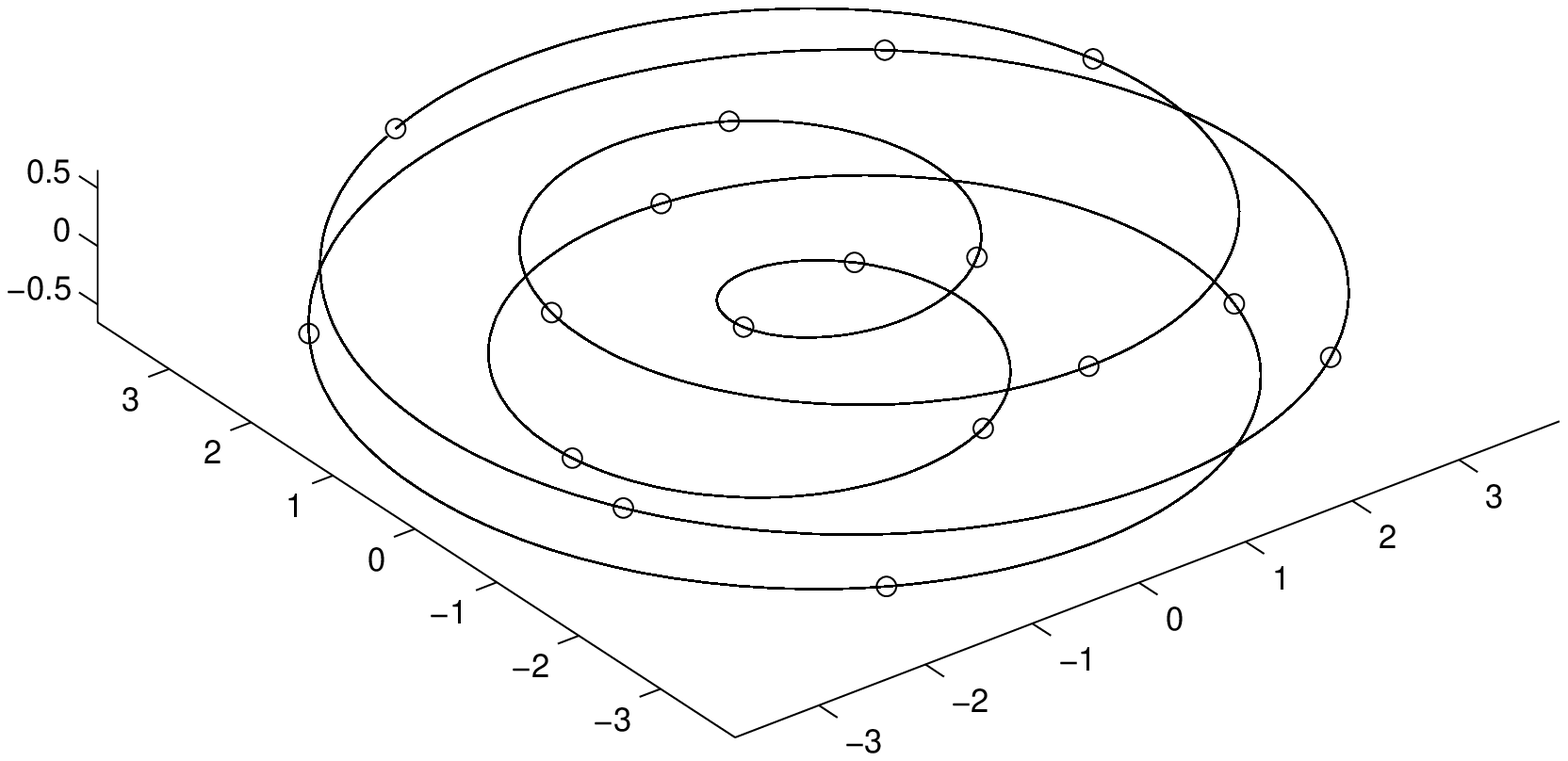,width=5.0cm}} \quad & \quad
  {\psfig{figure=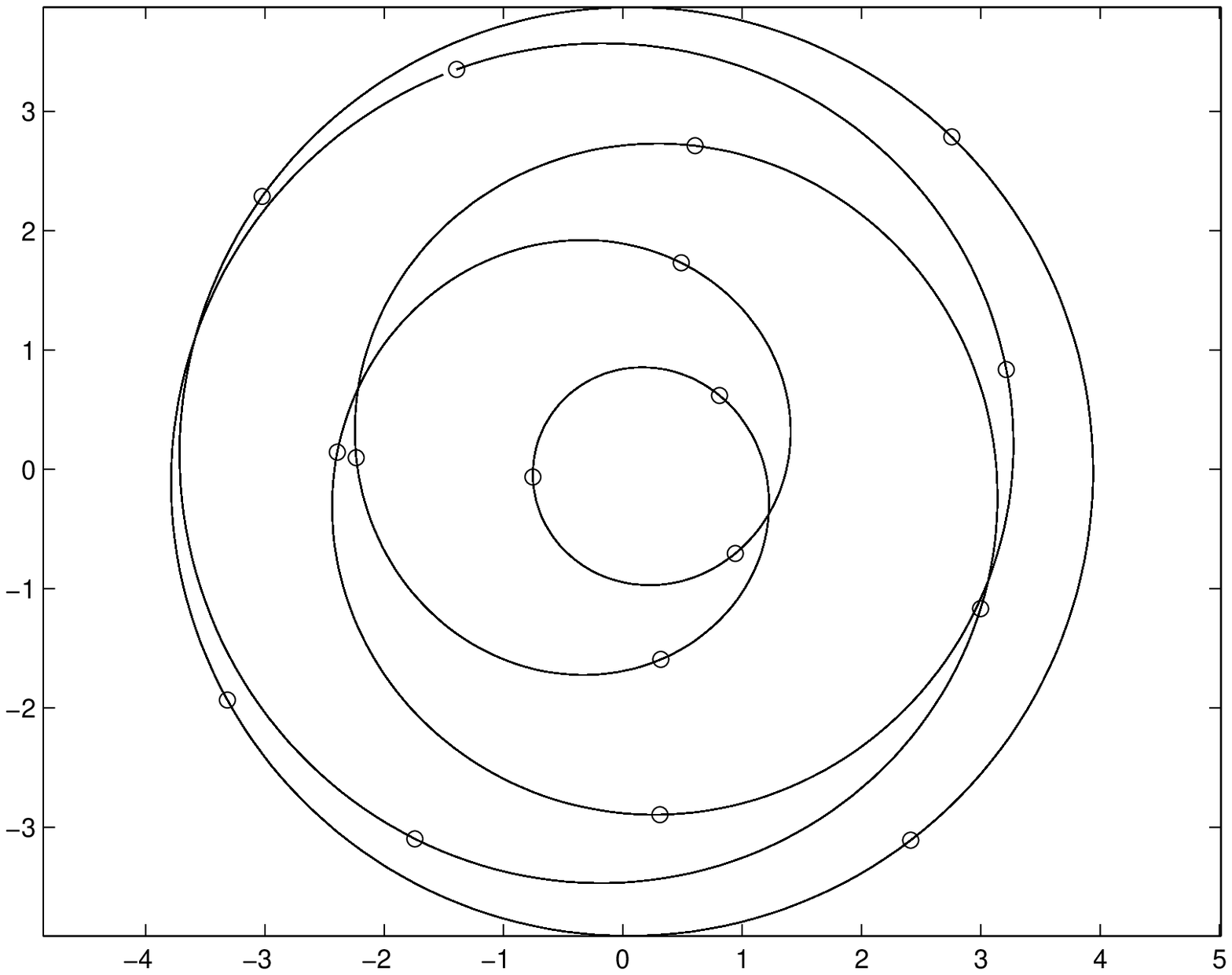,width=5.0cm}} \\
  {$n=17$, $\omega = 4.48$, $\alpha=1$} &
  \end{tabular} 
  \caption{Minima of the action with angular velocities close to the half on an integer.}
\end{center}
\end{figure}

The analisys of the minimizers of the action in a rotating system suggest the existence of saddle critical points; in fact in many cases we have an infinity of absolute minima, in other we proved the existence of absolute minima, but we can compute numerically  other relative minima. By a Mountain Pass algorithm we have found some example of such saddle critical points. In Figure 4 we show a saddle point for the 3-body problem with angular velocity $\omega = 1.5$, between the circles with minimal period $2\pi$ and $\pi$ (see Remark \ref{3corpi}). 
The second picture is a saddle point for the 4-body problem with angular velocity $\omega = 2.2$, between the circle with minimal period $2\pi/3$ (that is a local minimum) and orbit described in Theorem \ref{kdivn} for $k=2$. The last picture of Figure 4 is a saddle point for the 5-body problem with $\omega = 2.5$, between the circles with minimal period $\pi$ and $2\pi/3$. Remark that these saddle points are collisionless.

\begin{figure}[h!]
\begin{center}
  \begin{tabular}{cc}
  {\psfig{figure=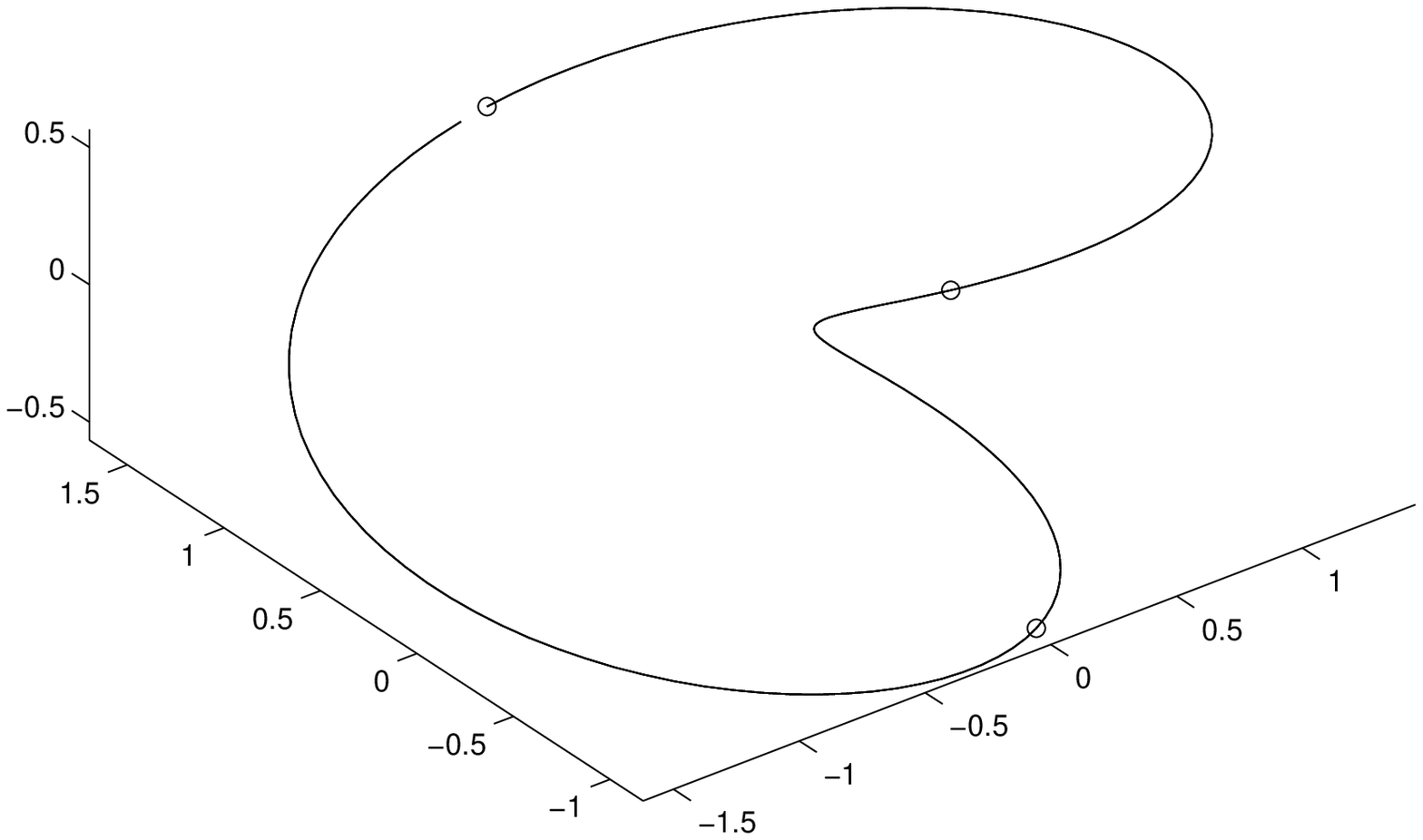,width=5.0cm}} \quad & \quad
  {\psfig{figure=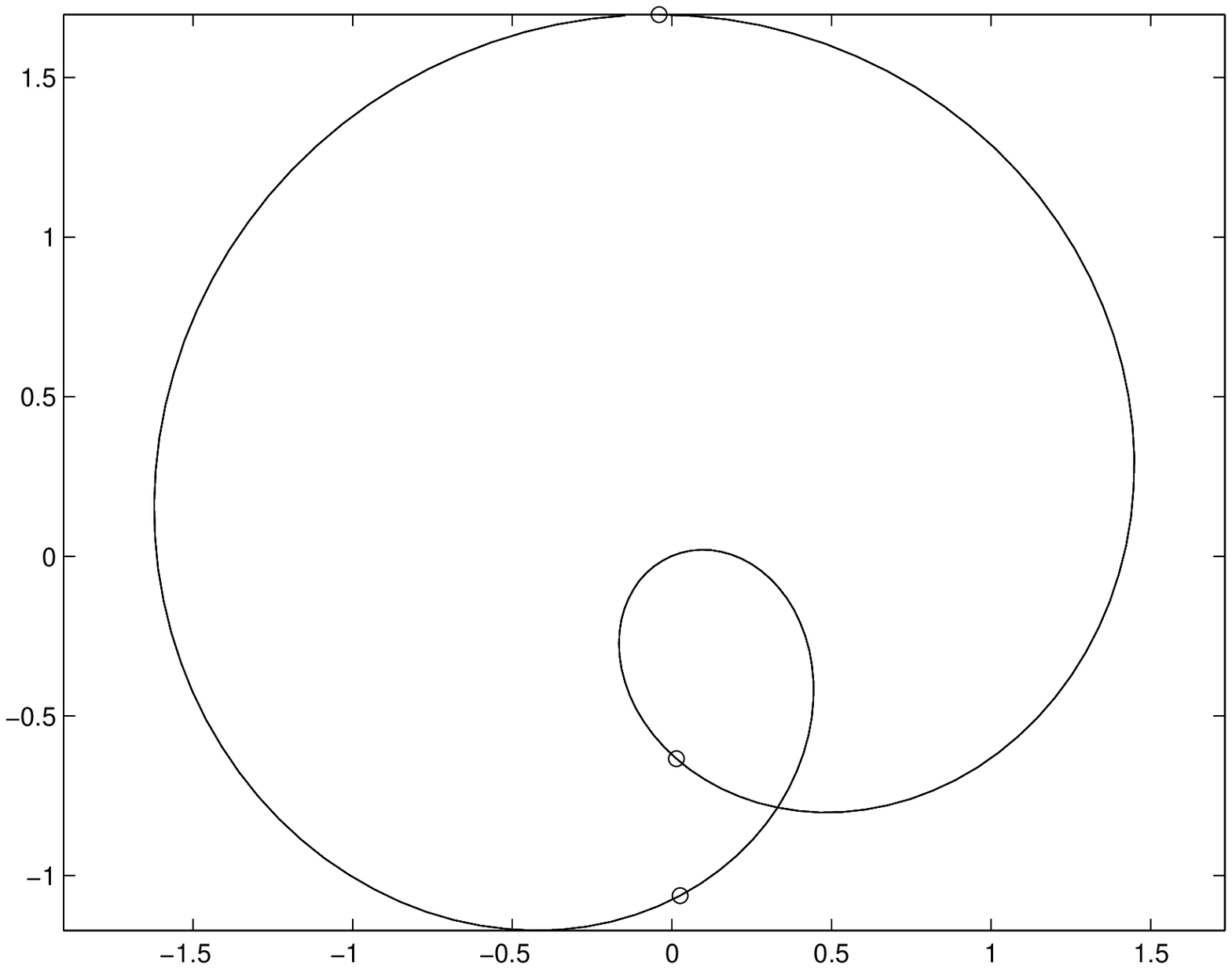,width=5.0cm}} \\
  {$n=3$, $\omega = 1.5$, $\alpha=1$} & \\
  {\psfig{figure=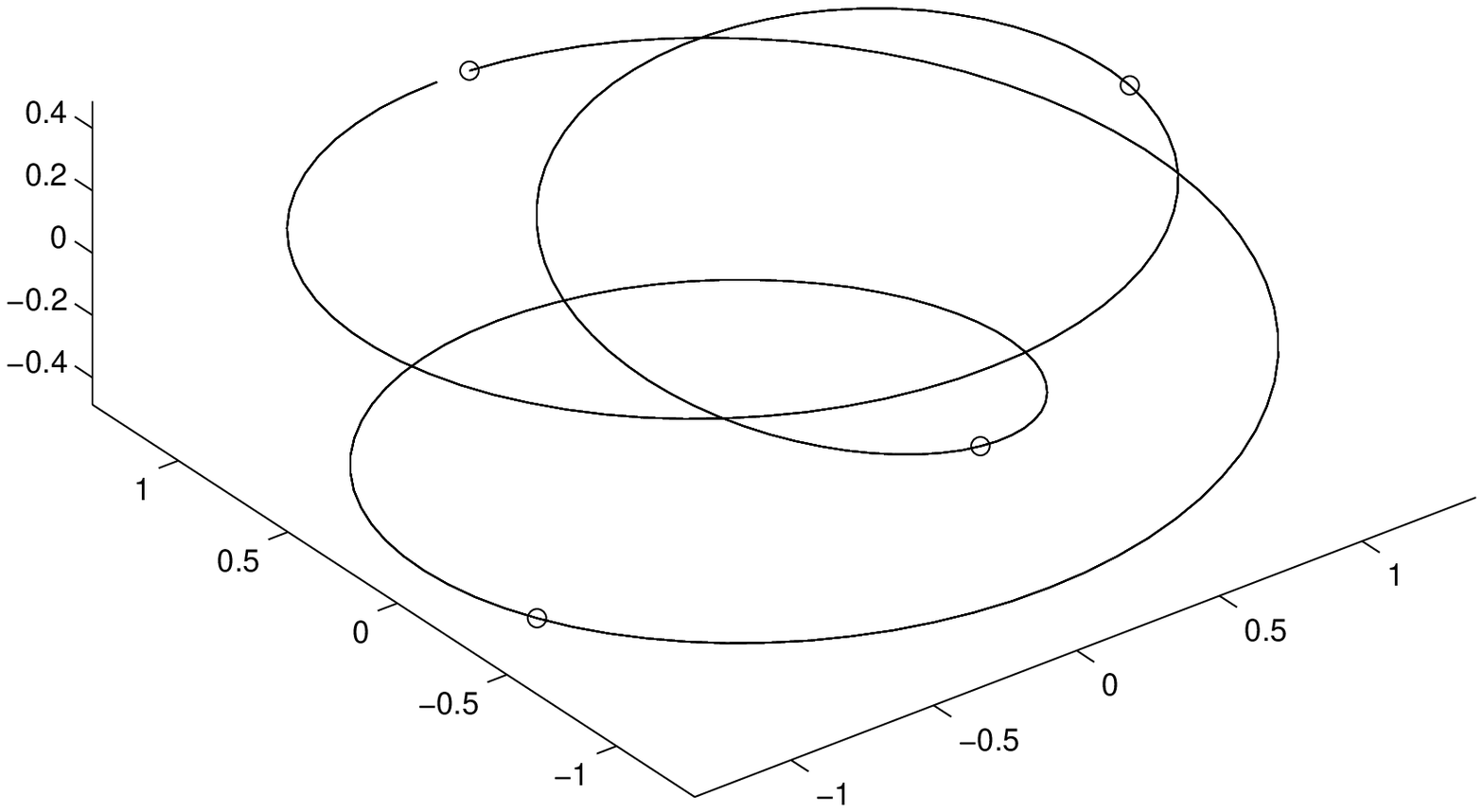,width=5.0cm}} \quad & \quad
  {\psfig{figure=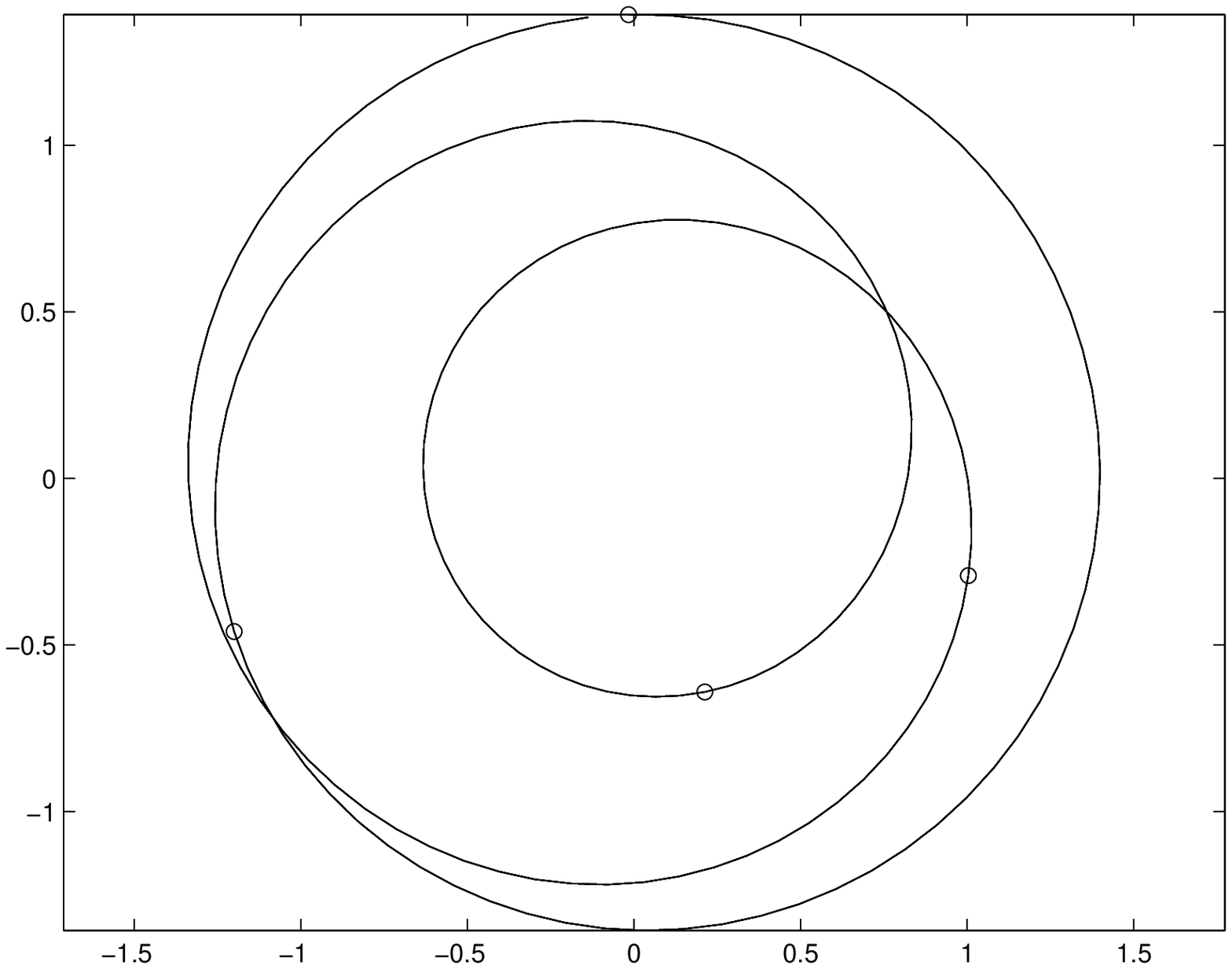,width=5.0cm}} \\
  {$n=4$, $\omega = 2.2$, $\alpha=1$} & \\
  {\psfig{figure=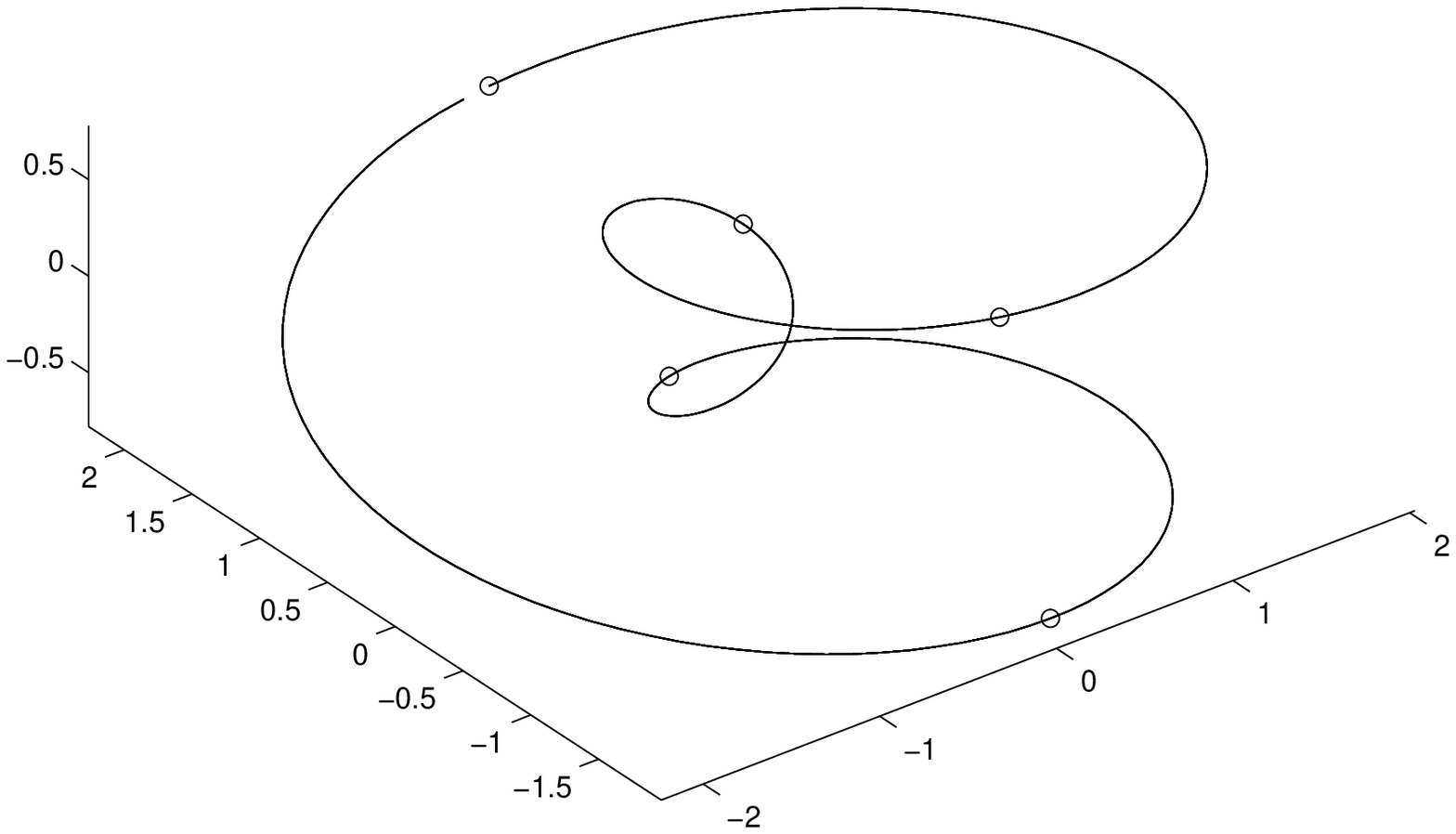,width=5.0cm}} \quad & \quad
  {\psfig{figure=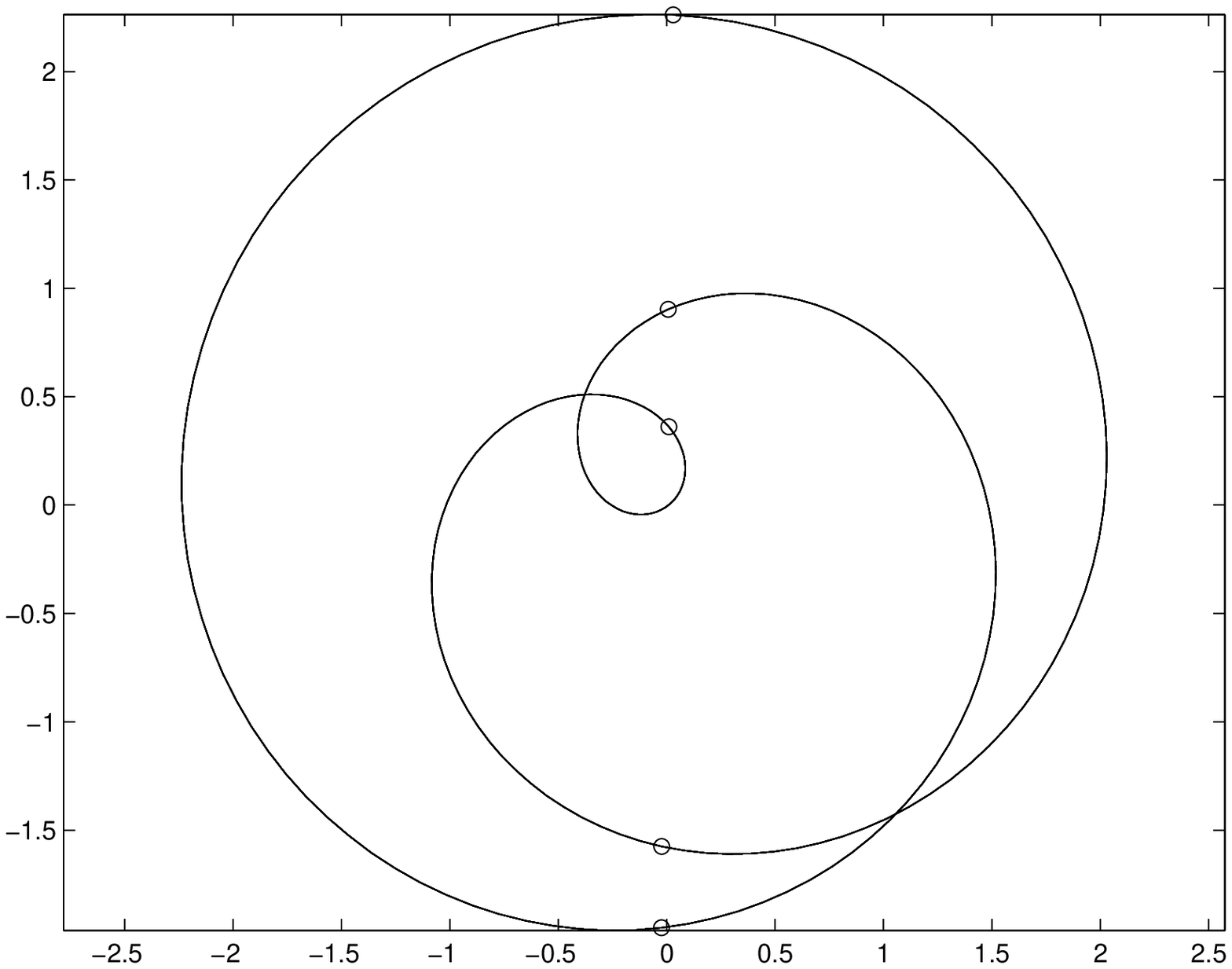,width=5.0cm}} \\
  {$n=5$, $\omega = 2.5$, $\alpha=1$} &  
  \end{tabular}
  \caption{Saddle points.}
\end{center}
\end{figure}

In this setting, we suggest a new way to interpretate the eight shaped curve that Chenciner and Montgomery found  in \cite{CM} as a minimum of the function integral among simple choreographies with some special symmetry constraints.\\
Consider $n$ bodies, $n$ odd, moving in $\RR^3$ on a loop $x(t)=(x_1(t),x_2(t),x_3(t))$ in the set $\Lambda$. We impose the following symmetry constraints
\begin{eqnarray}
\label{cs}
\left .
\begin{array}{ll}
	x_1(t) = x_1(t + \pi), & x_1(t) = -x_1(-t),\\
	x_2(t) = x_2(t + 2\pi), & x_2(t) = -x_2(-t),\\
  x_3(t) = x_3(t + 2\pi), & x_3(t) = x_3(-t).
\end{array}
 \right. 
\end{eqnarray}
Remark that these constraints are satisfied by a eight in the plane $x_3=0$ or by a circle in the plane $x_1=0$. In the set $\Lambda$ with (\ref{cs}), we have two circles that attains the minimun
$$  
(0,R\sin t,R\cos t) \quad \mbox{and} \quad (0,R\sin t,-R\cos t)
$$
where $R$ is determined in (\ref{radius}). This two curves are isolated minima; the Mountain Pass Theorem, ensures the existence of a saddle point, which is the minimum among the maximal points of the paths joining these two minima. Using a numerical mountain pass algorithm, we find that the eight in the plane $x_3=0$ is a saddle points between the two circles.

\end{document}